\documentclass[doc]{apa6}

\usepackage{amsmath}
\usepackage{amsfonts}
\usepackage{amssymb}
\usepackage{listings}
\makeatletter
\def\Let@{\def\\{\notag\math@cr}}
\makeatother
\usepackage{graphicx}
\usepackage[utf8]{inputenc}
\usepackage[T1]{fontenc}

\usepackage[utf8]{inputenc}
\usepackage{amsmath,amssymb}
\newcommand{\textVerb}[1]{\texttt{\mbox{#1}}}
\usepackage{subcaption}
\usepackage{algorithm,algorithmic}
\usepackage{color}

\usepackage[american]{babel}
\usepackage{csquotes}
\usepackage{natbib}
\usepackage{caption}
\usepackage{{url}}
\usepackage{rotating}

\newcommand\blfootnote[1]{%
  \begingroup
  \renewcommand\thefootnote{}\footnote{#1}%
  \addtocounter{footnote}{-1}%
  \endgroup
}

\title{Generalized Network Psychometrics: Combining Network and Latent Variable Models}

\shorttitle{GENERALIZED NETWORK PSYCHOMETRICS}

\author{Sacha Epskamp,\textsuperscript{1} Mijke T.\ Rhemtulla,\textsuperscript{2} and Denny Borsboom\textsuperscript{1}}

\twoaffiliations{1. University of Amsterdam: Department of Psychological Methods\\
	2. University of California, Davis: Department of Psychology}

\abstract{
We introduce the network model as a formal psychometric model, conceptualizing the covariance between psychometric indicators as resulting from pairwise interactions between observable variables in a network structure. This contrasts with standard psychometric models, in which the covariance between test items arises from the influence of one or more common latent variables. Here, we present two generalizations of the network model that encompass latent variable structures, establishing network modeling as parts of the more general framework of Structural Equation Modeling (SEM). In the first generalization, we model the covariance structure of latent variables as a network. We term this framework \emph{Latent Network Modeling} (LNM) and show that, with LNM, a unique structure of conditional independence relationships between latent variables can be obtained in an explorative manner. In the second generalization, the residual variance-covariance structure of indicators is modeled as a network. We term this generalization \emph{Residual Network Modeling} (RNM) and show that, within this framework, identifiable models can be obtained in which local independence is structurally violated. These generalizations allow for a general modeling framework that can be used to fit, and compare, SEM models, network models, and the RNM and LNM generalizations. This methodology has been implemented in the free-to-use software package \emph{lvnet}, which contains confirmatory model testing as well as two exploratory search algorithms: stepwise search algorithms for low-dimensional datasets and penalized maximum likelihood estimation for larger datasets. We show in simulation studies that these search algorithms performs adequately in identifying the structure of the relevant residual or latent networks. We further demonstrate the utility of these generalizations in an empirical example on a personality inventory dataset.
}

\begin{document}

\maketitle

\raggedbottom
\urlstyle{same}

\blfootnote{Manuscript accepted for publication in \emph{Psychometrika}.}

\section{Introduction}

Recent years have seen an emergence of network modeling in psychometrics \citep{borsboom2008psychometric, schmittmann2013}, with applications in clinical psychology (e.g., \citealt{van2015association,mcnally2015mental,fried2015}), psychiatry (e.g., \citealt{isvoranu,isvoranua}), health sciences (e.g., \citealt{kossakowski2015}), social psychology (e.g., \citealt{dalege2016toward,cramer2012dimensions}), and other fields (see for a review of recent literature \citealt{Fried_van_Borkulo_2016}). This line of literature stems from the network perspective of psychology \citep{cramer2010comorbidity},
in which noticeable macroscopic behavior---the co-occurrence of aspects of psychology such as cognitive abilities, psychopathological symptoms, or a set of behaviors---is hypothesized to not be due to the influence of unobserved common causes, such as general intelligence, psychopathological disorders, or personality traits, but rather to emergent behavior in a network of interacting psychological, sociological, and biological components. 
Network models are used to gain insight into this potentially high-dimensional interplay \citep{dynamics}. In practice, network models can be used as a sparse representation of the joint distribution of observed indicators, and as such these models show great promise in psychometrics by providing a perspective that complements latent variable modeling. Network modeling highlights variance that is unique to pairs of variables, whereas latent variable modeling focuses on variance that is shared across all variables \citep{costantini2015state}. As a result, network modeling and latent variable modeling can complement---rather than exclude---one-another. 

In this paper, we introduce the reader to this field of \emph{network psychometrics} \citep{netpsych} and formalize the network model for multivariate normal data, the Gaussian Graphical Model (GGM;  \citealt{lauritzen1996graphical}), as a formal psychometric model. We contrast the GGM to the Structural Equation Model (SEM;  \citealt{wright1921correlation,kaplan2008structural}) and show 
that the GGM can be seen as another way to approach modeling covariance structures as is typically done in psychometrics. In particular, rather than modeling the covariance matrix, the GGM models the \emph{inverse} of a covariance matrix. The GGM and SEM are thus very closely related: every GGM model and every SEM model imply a constrained covariance structure. We make use of this relationship to show that, through a reparameterization of the SEM model, the GGM model can be obtained in two different ways: first, as a network structure that relates a number of latent variables to each other, and second, as a network between residuals that remain given a fitted latent variable model. As such, the GGM can be modeled and estimated in SEM, which allows for network modeling of psychometric data to be carried out in a framework familiar to psychometricians and methodologists. In addition, this allows for one to assess the fit of a GGM, compare GGMs to one-another and compare a GGM to a SEM model.

However, the combination of GGM and SEM allows for more than fitting network models. As we will show, the strength of one framework can help overcome shortcomings of the other framework. In particular, SEM falls short in that exploratory estimation is complicated and there is a strong reliance on local independence, whereas the GGM falls short in that it assumes no latent variables. In this paper, we introduce network models for latent covariances and for residual covariances as two distinct generalized frameworks of both the SEM and GGM. The first framework, Latent Network Modeling (LNM), formulates a network among latent variables. This framework allows researchers to exploratively estimate conditional independence relationships between latent variables through model search algorithms; this estimation is difficult in the SEM framework due to the presence of equivalent models \citep{maccallum1993problem}. The second framework, which we denote Residual Network Modeling (RNM), formulates a network structure on the residuals of a SEM model. With this framework, researchers can circumvent critical assumptions of both SEM and the GGM: SEM typically relies on the assumption of local independence, whereas network modeling typically relies on the assumption that the covariance structure among a set of the items is not due to latent variables at all. The RNM framework allows researchers to estimate SEM models without the assumption of local independence (all residuals can be correlated, albeit due to a constrained structure on the inverse residual covariance matrix) as well as to estimate a network structure, while taking into account the fact that the covariance between items may be partly due to latent factors. 

While the powerful combination of SEM and GGM allows for confirmative testing of network structures both with and without latent variables, we recognize that few researchers have yet formulated strict confirmatory hypotheses in the relatively new field of network psychometrics. Often, researchers are more interested in exploratively searching a plausible network structure. To this end, we present two exploratory search algorithms. The first is a step-wise model search algorithm that adds and removes edges of a network as long as fit is improved, and the second uses penalized maximum likelihood estimation \citep{tibshirani1996regression} to estimate a sparse model. We evaluate the performance of these search methods in four simulation studies. Finally, the proposed methods have been implemented in a free-to-use R package, \textVerb{lvnet}, which we illustrate in an empirical example on personality inventory items \citep{psych}.

\section{Modeling Multivariate Gaussian Data}

Let  $\pmb{y}$ be the response vector of a random subject on $P$ items\footnote{Throughout this paper, vectors will be represented with lowercase boldfaced letters and matrices will be denoted by capital boldfaced letters. Roman letters will be used to denote observed variables and parameters (such as the number of nodes) and Greek letters will be used to denote latent variables and parameters that need to be estimated. The subscript $i$ will be used to denote the realized response vector of subject $i$ and omission of this subscript will be used to denote the response of a random subject.}. We assume  $\pmb{y}$ is centered and follows a multivariate Gaussian density:
\[
\pmb{y} \sim N_P\left( \pmb{0}, \pmb{\Sigma} \right),
\]
In which $\pmb{\Sigma}$ is a $P \times P$ variance--covariance matrix, estimated by some model-implied $\hat{\pmb{\Sigma}}$. Estimating $\hat{\pmb{\Sigma}}$ is often done through some form of \emph{maximum likelihood estimation}. If we measure $N$ independent samples of $\pmb{y}$ we can formulate the $N \times P$ matrix $\pmb{Y}$ containing realization $\pmb{y}^\top_i$ as its $i$th row. Let $\pmb{S}$ represent the sample variance--covariance matrix of $\pmb{Y}$:
\[
\pmb{S} = \frac{1}{N-1} \pmb{Y}^\top \pmb{Y}.
\]
In maximum likelihood estimation, we use $\pmb{S}$ to compute and minimize $-2$ times the log-likelihood function to find $\hat{\pmb{\Sigma}}$ \citep{lawley1940vi,joreskog1967general,regsem}:
\begin{equation}
\label{psychometrika:loglik}
\min_{\hat{\pmb{\Sigma}}} \left[ \log \det \left( \hat{\pmb{\Sigma}}\right)   + \mathrm{Trace}\left( \pmb{S} \hat{\pmb{\Sigma}}^{-1} \right)  - \log \det \left( \hat{\pmb{S}}\right)  - P \right].
\end{equation}
To optimize this expression, $\hat{\pmb{\Sigma}}$ should be estimated as closely as possible to $\pmb{S}$ and perfect fit is obtained if $\hat{\pmb{\Sigma}} = \pmb{S}$. A properly identified model with the same number of parameters ($K$) used to form $\hat{\pmb{\Sigma}}$ as there are unique elements in $\pmb{S}$ ($P(P+1)/2$ parameters) will lead to $\hat{\pmb{\Sigma}} = \pmb{S}$ and therefore a saturated model. The goal of modeling multivariate Gaussian data is to obtain some model for $\hat{\pmb{\Sigma}}$ with \emph{positive degrees of freedom}, $K < P(P+1)/2$, in which $\hat{\pmb{\Sigma}}$ resembles $\pmb{S}$ closely. 

\subsection{Structural Equation Modeling}

In Confirmatory Factor Analysis (CFA), $\pmb{Y}$ is typically assumed to be a causal linear effect of a set of $M$ centered latent variables, $\pmb{\eta}$, and independent residuals or error, $\pmb{\varepsilon}$:
\[
\pmb{y} = \pmb{\Lambda}\pmb{\eta} + \pmb{\varepsilon}.
\]
Here, $\pmb{\Lambda}$ represents a $P \times M$ matrix of \emph{factor loadings}. This model implies the following model for $\hat{\pmb{\Sigma}}$:
\begin{equation}
\label{psychometrika:eq:cfa}
\hat{\pmb{\Sigma}} = \pmb{\Lambda} \pmb{\Psi} \pmb{\Lambda}^\top + \pmb{\Theta},
\end{equation}
in which $\pmb{\Psi} = \mathrm{Var}\left( \pmb{\eta} \right)$ and $\pmb{\Theta} = \mathrm{Var}\left(\pmb{\varepsilon} \right)$. In Structural Equation Modeling (SEM), $\mathrm{Var}\left( \pmb{\eta} \right)$ can further be modeled by adding structural linear relations between the latent variables\footnote{We make use here of the convenient all-$y$ notation and do not distinguish between exogenous and endogenous latent variables \citep{hayduk1987structural}.}:
\[
\pmb{\eta} = \pmb{B} \pmb{\eta} + \pmb{\zeta},
\]
in which $\pmb{\zeta}$ is a vector of residuals and $\pmb{B}$ is an $M \times M$ matrix of regression coefficients. Now, $\hat{\pmb{\Sigma}}$ can be more extensively modeled as:
\begin{equation}
\label{psychometrika:eq:sem}
\hat{\pmb{\Sigma}} = \pmb{\Lambda}  \left( \pmb{I} - \pmb{B}   \right)^{-1} \pmb{\Psi} \left( \pmb{I} - \pmb{B}  \right)^{-1\top} \pmb{\Lambda}^{\top} + \pmb{\Theta},
\end{equation}
in which now $\pmb{\Psi} = \mathrm{Var}\left( \pmb{\zeta} \right)$. This framework can be used to model direct causal effects between observed variables by setting $\pmb{\Lambda} = \pmb{I}$ and $\pmb{\Theta} = \pmb{O}$, which is often called path analysis \citep{wright1934method}.

The $\pmb{\Theta}$ matrix is, like $\hat{\pmb{\Sigma}}$ and $\pmb{S}$, a $P \times P$ matrix; if $\pmb{\Theta}$ is fully estimated---contains no restricted elements---then $\pmb{\Theta}$ alone constitutes a saturated model. Therefore, to make either \eqref{psychometrika:eq:cfa} or \eqref{psychometrika:eq:sem} identifiable, $\pmb{\Theta}$ must be strongly restricted. Typically, $\pmb{\Theta}$ is set to be diagonal, a restriction often termed \emph{local independence} \citep{lord1968statistical, holland1986conditional} because indicators are independent of each other after conditioning on the set of latent variables. To improve fit, select off-diagonal elements of $\pmb{\Theta}$ can be estimated, but systematic violations of local independence---many nonzero elements in $\pmb{\Theta}$---are not possible as that will quickly make  \eqref{psychometrika:eq:cfa} and \eqref{psychometrika:eq:sem}  saturated or even over-identified. More precisely, $\pmb{\Theta}$ can \emph{not} be fully-populated---some elements of $\pmb{\Theta}$ must be set to equal zero---when latent variables are used. An element of $\pmb{\Theta}$ being fixed to zero indicates that two variables are \emph{locally independent} after conditioning on the set of latent variables.  As such, local independence is a \emph{critical assumption} in both CFA and SEM; if local independence is systematically violated, CFA and SEM will never result in correct models.

The assumption of local independence has led to critiques of the factor model and its usage in psychology; local independence appears to be frequently violated due to direct causal effects, semantic overlap, or reciprocal interactions between putative indicators of a latent variable \citep{borsboom2008psychometric, cramer2010comorbidity,borsboom2011small, cramer2012dimensions, schmittmann2013}. In psychopathology research, local independence of symptoms given a person's level of a latent mental disorder has been questioned \citep{borsboom2013network}. For example, three problems associated with depression are ``fatigue'', ``concentration problems'' and ``rumination''. It is plausible that a person who suffers from fatigue will also concentrate more poorly, as a direct result of being fatigued and regardless of his or her level of depression. Similarly, rumination might lead to poor concentration. In another example, \citet{kossakowski2015} describe the often-used SF-36 questionnaire \citep{ware1992} designed to measure health related quality of life. The SF-36 contains items such as ``can you walk for more than one kilometer'' and ``can you walk a few hundred meters''. Clearly, these items can never be locally independent after conditioning on any latent trait, as one item (the ability to walk a few hundred meters) is a \emph{prerequisite} for the other (walking more than a kilometer). In typical applications, the excessive covariance between items of this type is typically left unmodeled, and treated instead by combining items into a subscale or total score that is subsequently subjected to factor analysis; of course, however, this is tantamount to ignoring the relevant psychometric problem rather than solving it.

Given the many theoretically expected violations of local independence in psychometric applications, many elements of $\pmb{\Theta}$ in both \eqref{psychometrika:eq:cfa} and \eqref{psychometrika:eq:sem} should ordinarily be freely estimated. Especially when violations of local independence are expected to be due to causal effects of partial overlap, residual correlations should not be constrained to zero; in addition, a chain of causal relationships between indicators can lead to all residuals to become correlated. Thus, even when latent factors cause much of the covariation between measured items, fitting a latent variable model that involves local independence may not fully account for correlation structure between measured items. Of course, in practice, many psychometricians are aware of this problem, which is typically addressed by freeing up correlations between residuals to improve model fit. However, this is usually done in an ad-hoc fashion, on the basis of inspection of modification indices and freeing up error covariances one by one, which is post hoc, suboptimal, and involves an uncontrolled journey through the model space. As a result, it is often difficult to impossible to tell how exactly authors arrived at their final reported models. As we will show later in this paper, this process can be optimized and systematized using network models to connect residuals on top of a latent variable structure.

\subsection{Network Modeling}

Recent authors have suggested that the potential presence of causal relationships between measured variables may allow the explanation of the covariance structure without the need to invoke any latent variables \citep{borsboom2008psychometric, cramer2010comorbidity,borsboom2011small, schmittmann2013}. The interactions between indicators can instead be modeled as a \emph{network}, in which indicators are represented as nodes that are connected by edges representing \emph{pairwise interactions}. Such interactions indicate the presence of covariances that cannot be explained by any other variable in the model and can represent---possibly reciprocal---causal relationships. Estimating a network structure on psychometric data is termed \emph{network psychometrics} \citep{netpsych}. Such a network of interacting components can generate data that fit factor models well, as is commonly the case in psychology. \citet{van2006dynamical} showed that the positive manifold of intelligence---which is commonly explained with the general factor for intelligence, $g$---can emerge from a network of mutually benefiting cognitive abilities. \citet{borsboom2011small} showed that a network of psychopathological symptoms, in which disorders are modeled as clusters of symptoms, could explain comorbidity between disorders. Furthermore, \citet{netpsych} showed that the Ising model for ferromagnetism \citep{ising1925beitrag}, which models magnetism as a network of particles, is equivalent to multidimensional item response theory \citep{reckase2009multidimensional}.

\begin{figure}
\centering
\includegraphics[width=0.5\textwidth]{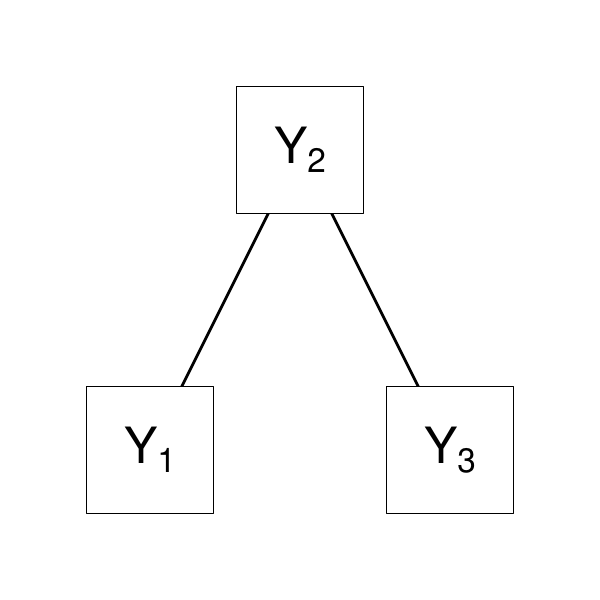}
\caption{Example of a pairwise Markov Random Field model. Edges in this model indicate pairwise interactions, and are drawn using undirected edges to distinguish from (bidirectional)  covariances. Rather than a model for  marginal  associations  (such  as a network indicating covariances), this is a model for \emph{conditional} associations.  The network above encodes that $Y_1$ and $Y_3$ are independent after conditioning on $Y_2$. Such a model allows all three variables to correlate while retaining one degree of freedom (the model only has two parameters)}
\label{psychometrika:fig:UN}
\end{figure}

In network psychometrics, psychometric data are modeled through directed or undirected networks. Directed networks are equivalent to path analysis models. For modeling undirected networks, pairwise Markov Random Fields \citep{lauritzen1996graphical, murphy2012machine} are used. In these models, each variable is represented by a node, and nodes are connected by a set of edges. If two nodes, $y_j$ and $y_k$, are not connected by an edge, then this means they are independent after conditioning on the set of all other nodes, $\pmb{y}^{-(j,k)}$. Whenever two nodes cannot be rendered independent conditional on the other nodes in the system, they are said to feature in a \emph{pairwise interaction}, which is represented by an \emph{undirected} edge---an edge with no arrows---to contrast such an effect from covariances typically represented in the SEM literature with bidirectional edges. Figure~\ref{psychometrika:fig:UN} represents such a network model, in which nodes  $y_1$ and $y_3$ are independent after conditioning on node $y_2$. Such a model can readily arise from direct interactions between the nodes. For example, this conditional independence structure would emerge if $y_2$ is a common cause of $y_1$ and $y_3$, or if $y_2$ is the mediator in a causal path between $y_1$ and $y_3$. In general, it is important to note that pairwise interactions are not mere correlations; two variables may be strongly correlated but unconnected (e.g., when both are caused by another variable in the system) and they may be uncorrelated but strongly connected in the network (e.g., when they have a common effect in the system). For instance, in the present example the model does \emph{not} indicate that $y_1$ and $y_3$ are uncorrelated, but merely indicates that any correlation between $y_1$ and $y_3$ is due to their mutual interaction with $y_2$; a network model in which either directly or indirectly connected paths exist between all pairs of nodes typically implies a fully populated (no zero elements) variance--covariance matrix. 
 
In the case of multivariate Gaussian data, this model is termed the Gaussian Graphical Model (GGM; \citealt{lauritzen1996graphical}). Here, the \emph{partial correlation coefficient} is sufficient to test the degree of conditional independence of two variables after conditioning on all other variables; if the partial correlation coefficient is zero, there is conditional independence and hence no edge in the network. As such, partial correlation coefficients can directly be used in the network as \emph{edge weights}; the strength of connection between two nodes\footnote{A saturated GGM is also called a partial correlation network because it contains the sample partial correlation coefficients as edge weights.}. Such a network is typically encoded in a symmetrical and real valued $p \times p$ \emph{weight matrix}, $\pmb{\Omega}$, in which element $\omega_{jk}$ represents the edge weight between node $j$ and node $k$:
\[
\mathrm{Cor}\left(y_j ,y_k \mid \pmb{y}^{-(j,k)}\right) = \omega_{jk} = \omega_{kj}.
\]
The partial correlation coefficients can be directly obtained from the inverse of variance--covariance matrix $\hat{\pmb{\Sigma}}$, also termed the \emph{precision matrix} $\hat{\pmb{K}}$ \citep{lauritzen1996graphical}:
\[
\mathrm{Cor}\left( y_j, y_k \mid \pmb{y}^{-(j,k)}\right)  =- \frac{  \kappa_{jk} }{ \sqrt{\kappa_{kk}}\sqrt{\kappa_{jj}} }.
\]
Thus, element $\kappa_{jk}$ of the precision matrix is proportional to to the partial correlation coefficient of variables $y_j$ and $y_k$ after conditioning on all other variables. Since this process simply involves standardizing the precision matrix, we propose the following model\footnote{To our knowledge, the GGM has not yet been framed in this form. We chose this form because it allows for clear modeling and interpretation of the network parameters.}:
\begin{equation}
\label{psychometrika:netmodel}
\hat{\pmb{\Sigma}} = \hat{\pmb{K}}^{-1} = \pmb{\Delta} \left( \pmb{I} - \pmb{\Omega}  \right)^{-1}  \pmb{\Delta},
\end{equation}
in which $\pmb{\Delta}$ is a diagonal matrix with $\delta_{jj} =\kappa^{-\frac{1}{2}}_{jj}$ and $ \pmb{\Omega}$ has zeroes on the diagonal. This model allows for confirmative testing of the GGM structures on psychometric data. Furthermore, the model can be compared to a saturated model (fully populated off-diagonal values of $ \pmb{\Omega}$) and the independence model ($ \pmb{\Omega} = \pmb{O}$), allowing one to obtain $\chi^2$ fit statistics as well as fit indices such as the RMSEA \citep{browne1992alternative}  and CFI \citep{bentler1990comparative}. Such methods of assessing model fit have not yet been used in network psychometrics.

Similar to CFA and SEM, the GGM relies on a critical assumption; namely, that covariances between observed variables are \emph{not} caused by any latent or unobserved variable. If we estimate a GGM in a case where, in fact, a latent factor model was the true data generating structure, then generally we would expect the GGM to be saturated---i.e., there would be no missing edges in the GGM \citep{chandrasekaran2010latent}. A missing edge in the GGM indicates the presence of conditional independence between two indicators given all other indicators; we do not expect indicators to become independent given subsets of other indicators (see also \citealt{ellis1997tail}; \citealt{holland1986conditional}). Again, this critical assumption might not be plausible. While variables such as ``Am indifferent to the feelings of others'' and ``Inquire about others' well-being'' quite probably interact with each other, it might be far-fetched to assume that no unobserved variable, such as a personality trait, in part also causes some of the variance in responses on these items.

\section{Generalizing Factor Analysis and Network Modeling}

We propose two generalizations of both SEM and the GGM that both allow the modeling of network structures in SEM.  In the first generalization, we adopt the CFA\footnote{We use the CFA framework instead of the SEM framework here as the main application of this framework is in exploratively estimating relationships between latent variables.} decomposition in $\eqref{psychometrika:eq:cfa}$ and model the variance--covariance matrix of latent variables as a GGM:
\begin{align*}
\pmb{\Psi} &=  \pmb{\Delta}_{\pmb{\Psi}} \left( \pmb{I} - \pmb{\Omega}_{\pmb{\Psi}} \right)^{-1} \pmb{\Delta}_{\pmb{\Psi}}.
\end{align*}
This framework can be seen as modeling conditional independencies between latent variables not by directed effects (as in SEM) but as an undirected network. As such, we term this framework \emph{latent network modeling} (LNM). 

In the second generalization, we adopt the SEM decomposition of the variance--covariance matrix in \eqref{psychometrika:eq:sem} and allow the residual variance--covariance matrix $\pmb{\Theta}$ to be modeled as a GGM:
\begin{align*}
\pmb{\Theta} &=  \pmb{\Delta}_{\pmb{\Theta}} \left( \pmb{I} - \pmb{\Omega}_{\pmb{\Theta}} \right)^{-1} \pmb{\Delta}_{\pmb{\Theta}}.
\end{align*}
Because this framework conceptualizes associations between residuals as pairwise \emph{interactions}, rather than correlations, we term this framework \emph{Residual Network Modeling} (RNM). Using this framework allows---as will be described below---for a powerful way of fitting a confirmatory factor structure even though local independence is systematically violated and all residuals are correlated.

\begin{figure}
    \centering
    \begin{subfigure}[b]{0.45\textwidth}
        \includegraphics[width=\textwidth]{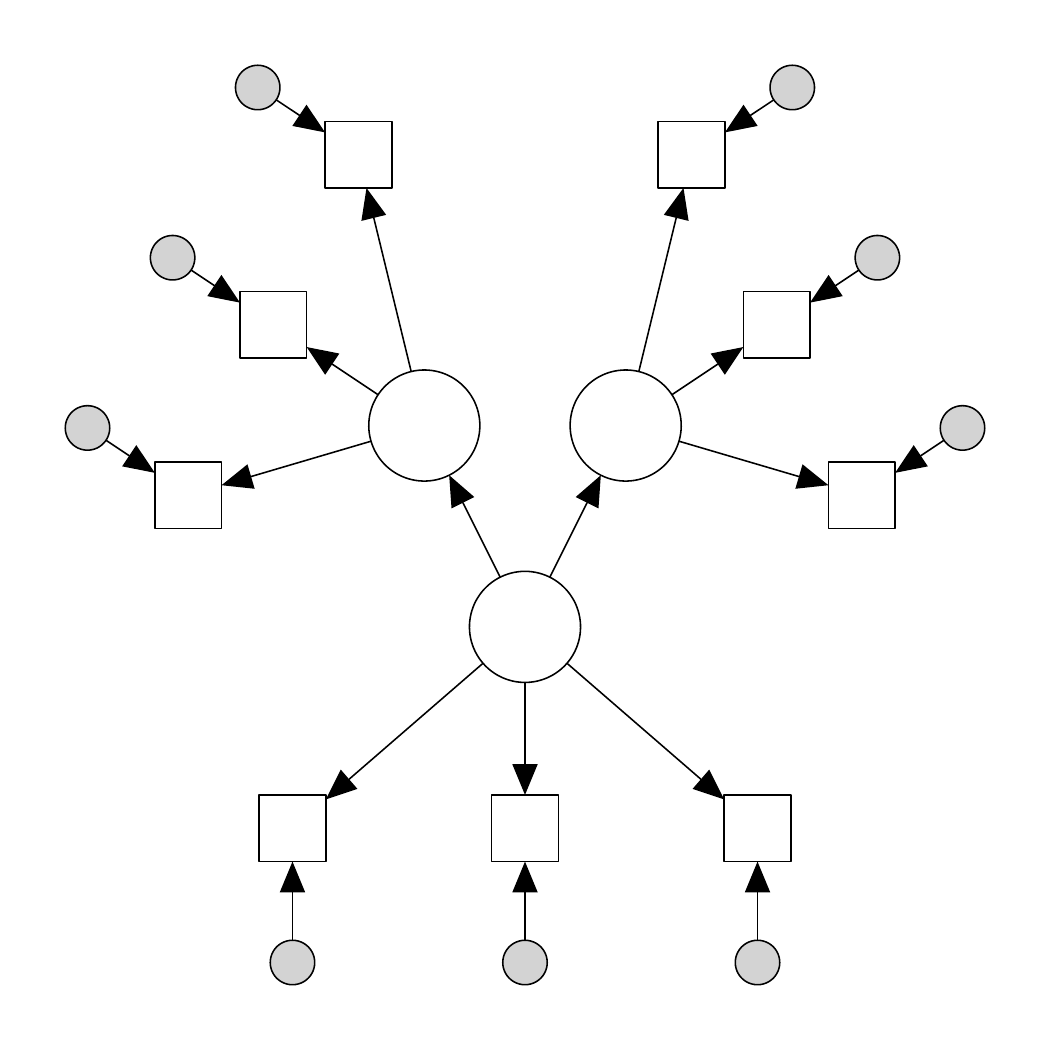}
        \caption{Structural Equation Modeling}
        \label{fig:sem}
    \end{subfigure}
    \begin{subfigure}[b]{0.45\textwidth}
        \includegraphics[width=\textwidth]{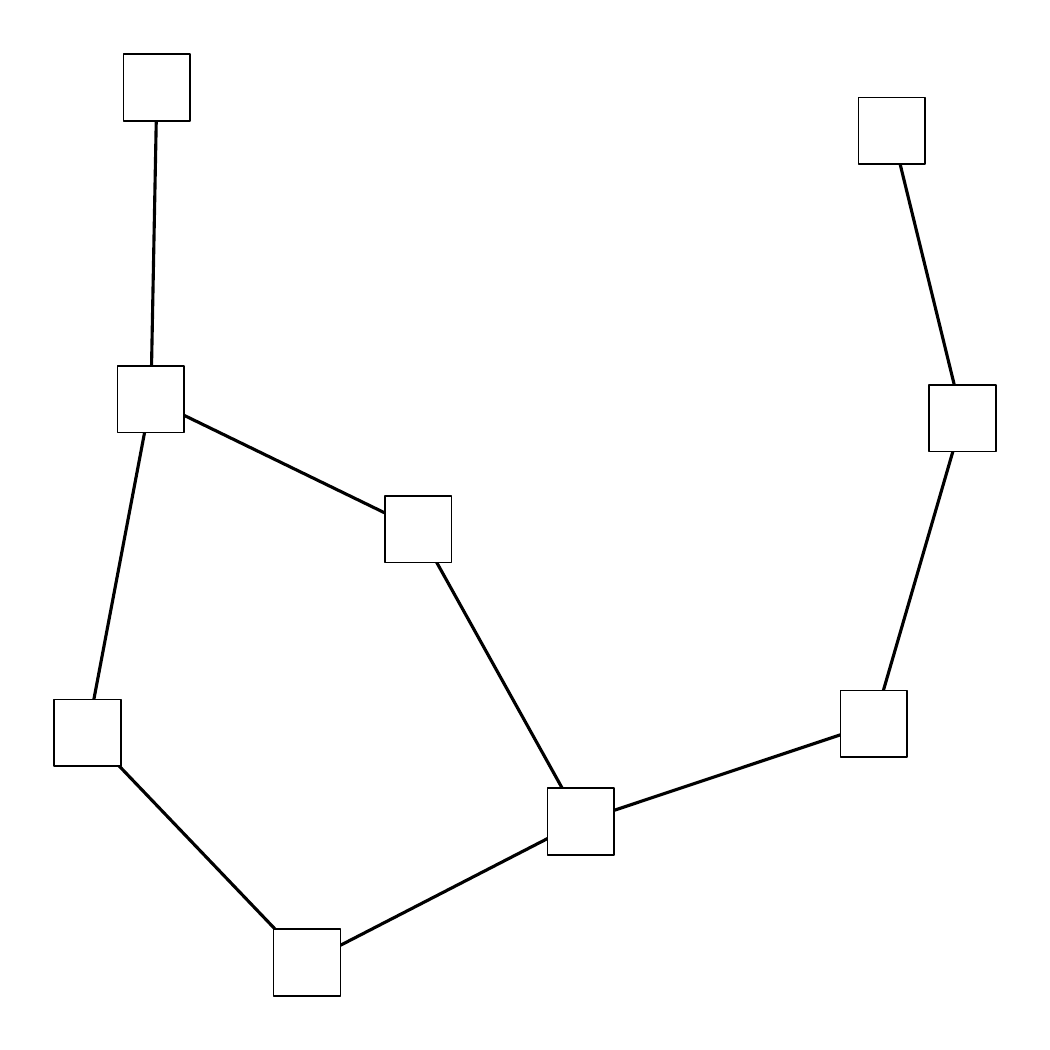}
        \caption{Network Modeling}
        \label{fig:NM}
    \end{subfigure} \\
    \begin{subfigure}[b]{0.45\textwidth}
        \includegraphics[width=\textwidth]{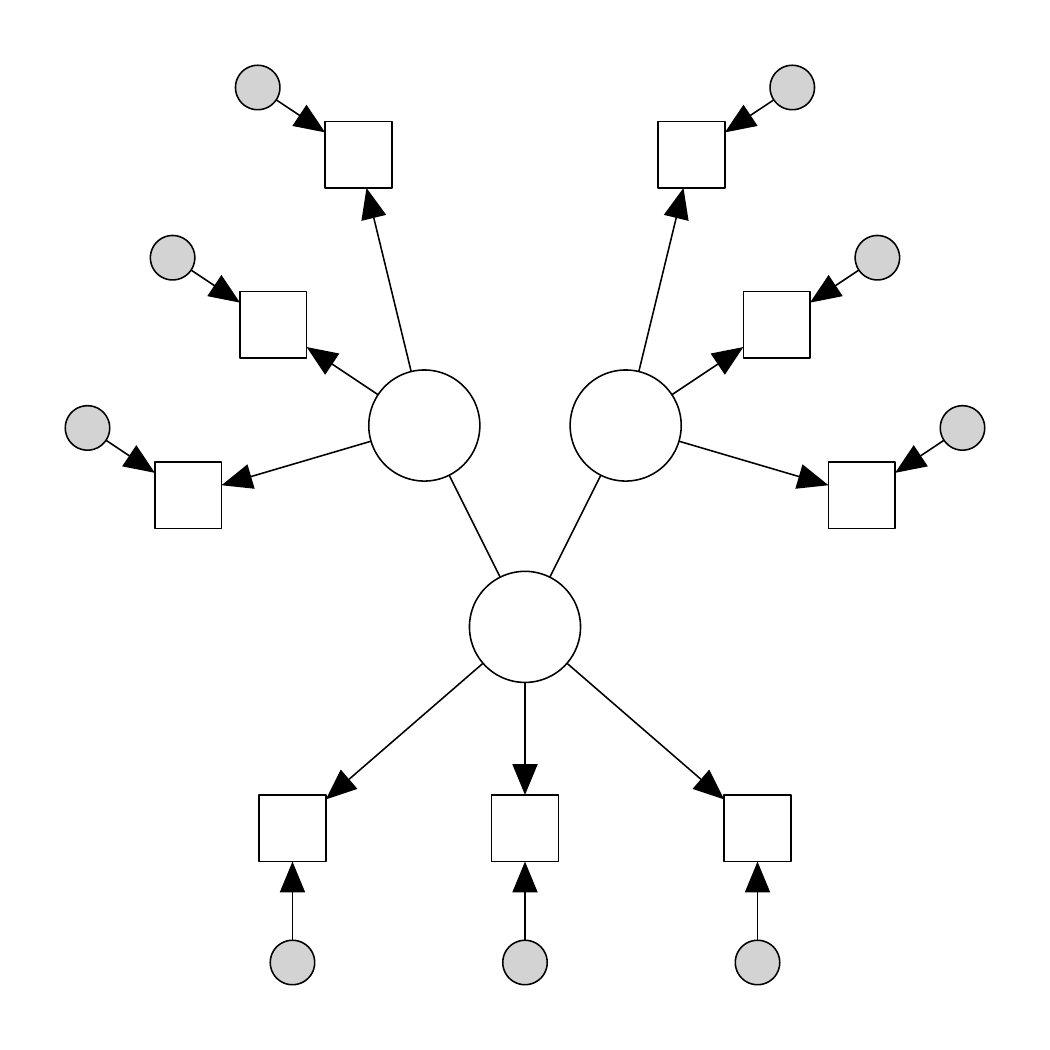}
        \caption{Latent Network Modeling}
        \label{fig:LNM}
    \end{subfigure}  
    \begin{subfigure}[b]{0.45\textwidth}
        \includegraphics[width=\textwidth]{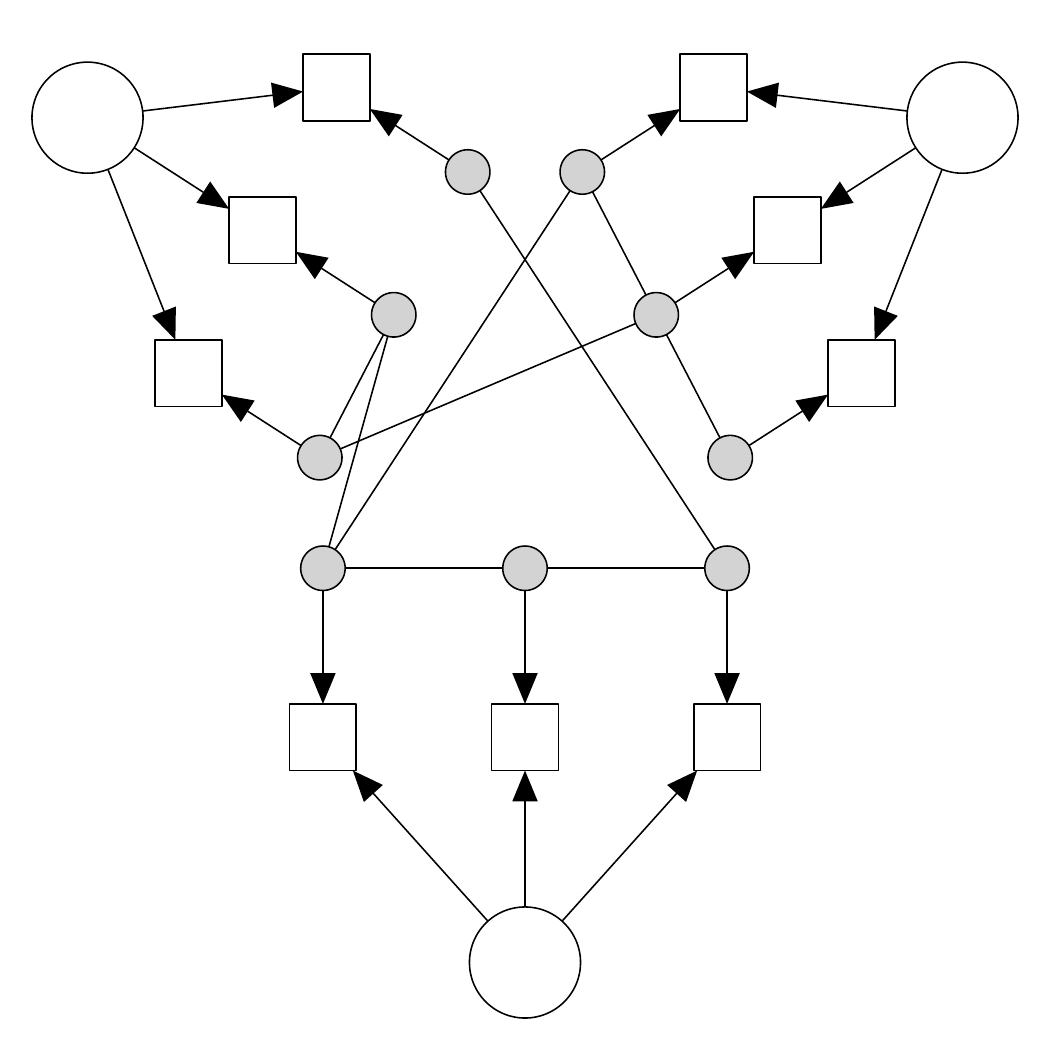}
        \caption{Residual Network Modeling}
        \label{fig:RNM}
    \end{subfigure} 
    \caption{Examples of possible models under four different modeling frameworks. Circular nodes indicate latent variables, square nodes indicate manifest variables and gray nodes indicate residuals. Directed edges indicate factor loadings or regression parameters and undirected edges indicate pairwise interactions. Note that such undirected edges do \emph{not} indicate covariances, which are typically denoted with bidirectional edges. Replacing covariances with interactions is where the network models differ from typical SEM.}\label{psychometrika:fig:models}
\end{figure}

Figure~\ref{psychometrika:fig:models} shows four different examples of possible models that are attainable under the SEM, LNM and RNM frameworks. Panel~(a) shows a typical SEM model in which one latent variable functions as a common cause of two others. Panel~(b) shows a network model which can be estimated using both the RNM and the LNM frameworks. Panel~(c) shows a completely equivalent LNM model to the SEM model of Panel~(a) in which the direction of effect between latent variables is not modeled. Finally, panel~(d) shows a model in which three exogenous latent variables underlie a set of indicators of which the residuals form a network. The remainder of this section will describe RNM and LNM in more detail and will outline the class of situations in which using these models is advantageous over CFA or SEM.

\subsection{Latent Network Modeling}

The LNM framework models the latent variance--covariance matrix of a CFA model as a GGM:
\begin{equation}
\label{psychometrika:eq:LNM}
\hat{\pmb{\Sigma}} = \pmb{\Lambda} \pmb{\Delta}_{\pmb{\Psi}} \left( \pmb{I} - \pmb{\Omega}_{\pmb{\Psi}} \right)^{-1} \pmb{\Delta}_{\pmb{\Psi}} \pmb{\Lambda}^\top + \pmb{\Theta}.
\end{equation}
This allows researchers to model conditional independence relationships between latent variables without making the implicit assumptions of directionality or acyclicness. In SEM, $\pmb{B}$ is typically modeled as a directed acyclic graph (DAG), meaning that elements of $\pmb{B}$ can be represented by directed edges, and following along the path of these edges it is not possible to return to any node (latent variable). The edges in such a DAG can be interpreted as causal, and in general they imply a specific set of conditional independence relationships between the nodes \citep{pearl2000causality}. 

While modeling conditional independence relationships between latent variables as a DAG is a powerful tool for testing strictly confirmatory hypotheses, it is less useful for more exploratory estimation. Though there have been recent advances in exploratory estimation of DAGs within an SEM framework (e.g., \citealt{gates2012group,rosa2012post}), many equivalent DAGs can imply the same conditional independence relationships, and thus fit the data equally well even though their causal interpretation can be strikingly different \citep{maccallum1993problem}. Furthermore, the assumption that the generating model is acyclic---which, in practice, often is made on purely pragmatic grounds to identify a model---is problematic in that much psychological behavior can be assumed to have at least some cyclic and complex behavior and feedback \citep{schmittmann2013}. Thus, the true conditional independence relationships in a dataset can lead to many equivalent compositions of $\pmb{B}$, and possibly none of them are the true model.

\begin{figure}
        \centering
        \begin{subfigure}{0.24\textwidth}
                \includegraphics[width=1\textwidth,page=1]{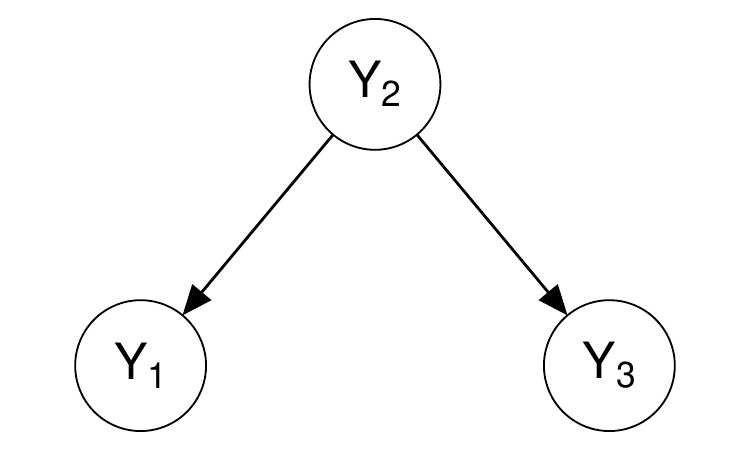}
                \caption{}
                \label{GGMDAG1}
        \end{subfigure}%
               \begin{subfigure}{0.24\textwidth}
                \includegraphics[width=1\textwidth,page=3]{Figure3.pdf}
                \caption{}
                \label{GGMDAG3}
        \end{subfigure}%
        \begin{subfigure}{0.24\textwidth}
                \includegraphics[width=1\textwidth,page=5]{Figure3.pdf}
                \caption{}
                \label{GGMDAG5} 
         \end{subfigure}
              \begin{subfigure}{0.24\textwidth}
                \includegraphics[width=1\textwidth,page=6]{Figure3.pdf}
                \caption{}
                \label{GGMDAG7} 
        \end{subfigure} \\
        \begin{subfigure}{0.24\textwidth}
                \includegraphics[width=1\textwidth,page=2]{Figure3.pdf}
                \caption{}
                \label{GGMDAG2}
        \end{subfigure}%
                \begin{subfigure}{0.24\textwidth}
                \includegraphics[width=1\textwidth,page=4]{Figure3.pdf}
                \caption{}
                \label{GGMDAG4}
        \end{subfigure}%
        \begin{subfigure}{0.24\textwidth}
                \includegraphics[width=1\textwidth,page=6]{Figure3.pdf}
                \caption{}
                \label{GGMDAG6} 
        \end{subfigure}        
        \begin{subfigure}{0.24\textwidth}
                \includegraphics[width=1\textwidth,page=8]{Figure3.pdf}
                \caption{}
                \label{GGMDAG8} 
        \end{subfigure}
        \caption{Equivalent models between directed acyclic graphs (DAG; top) and Gaussian Graphical Models (GGM; bottom). Each column of graphs show two models that are equivalent. Panels (a), (b), (e) and (f) all represent the same conditional independence structure: $Y_1$ and $Y_3$ are independent after conditioning on $Y_2$. Panel~(c) represents that  $Y_1$ and $Y_3$ are marginally independent even though they both are correlated with $Y_2$, a structure that cannot be represented in a GGM with only two edges. Panel (h) shows that $Y_1$ and $Y_3$ are independent after conditioning on the set $Y_2$ and $Y_4$, which cannot be represented with a DAG of four edges.}
        \label{psychometrika:fig:GGMvsDAG}
\end{figure}

In psychometrics and SEM the GGM representation has not been very prominent, even though it has some manifest benefits over the attempt to identify DAGs directly. For example, by modeling conditional independence relationships between latent variables as a GGM, many relationships can be modeled in a simpler way as compared to a DAG. In addition, in the GGM each set of conditional independence relations only corresponds to one model: there are no equivalent GGMs with the same nodes. Figure~\ref{psychometrika:fig:GGMvsDAG} shows a comparison of several conditional independence relations that can be modeled equivalently or not by using a GGM or by using a DAG. Panel~(a) and Panel~(b) show two DAGs that represent the same conditional independence relations, $y_1  \perp\!\!\!\perp y_3 \mid y_2$, which can both be represented by the same GGM shown in Panel~(e) and Panel~(f). There are some conditional independence relations that a GGM cannot represent in the same number of parameters as a DAG; Panel~(c) shows a collider structure that cannot be exactly represented by a GGM (the best fitting GGM would feature three edges instead of two). On the other hand, there are also conditional independence relationships that a GGM can represent and a DAG cannot; the cycle of Panel~(h) cannot be represented by a DAG. Further equivalences and differences between GGMs and DAGs are beyond the scope of this paper, but haven been well described in the literature (e.g., chapter 3 of \citealt{lauritzen1996graphical}; \citealt{koller2009probabilistic}; \citealt{kolaczyk2009statistical}). In sum, the GGM offers a natural middle ground between zero-order correlations and DAGs: every set of zero-order correlations implies exactly one GGM, and every DAG implies exactly one GGM. In a sense, the road from correlations to DAGs (including hierarchical factor models) thus always must pass through the realm of GGMs, which acts as a bridge between the correlational and causal worlds. 

Because there are no equivalent undirected models possible, LNM offers a powerful tool for exploratory estimation of relationships between latent variables. For example, suppose one encounters data generated by the SEM model in Figure~\ref{psychometrika:fig:models}, Panel~(a). Without prior theory on the relations between latent variables, exploratory estimation on this dataset would lead to three completely equivalent models: the one shown in Figure~\ref{psychometrika:fig:models}, Panel~(c) and two models in which the common cause instead is the middle node in a causal chain. As the number of latent variables increases, the potential number of equivalent models that encode the same conditional independence relationships grows without bound. The LNM model in Panel~(c) of Figure~\ref{psychometrika:fig:models} portrays the same conditional independence relationship as the SEM model in Panel~(a) of Figure~\ref{psychometrika:fig:models}, while having no equivalent model. Exploratory estimation could easily find this model, and portrays the retrieved relationship in a clear and unambiguous way.

A final benefit of using LNM models is that they allow network analysts to construct a network while taking measurement error into account. So far, networks have been constructed based on single indicators only and no attempt has been made to remediate measurement error. By forming a network on graspable small concepts measured by a few indicators, the LNM framework can be used to control for measurement error.

\subsection{Residual Network Modeling}

In the RNM framework the residual structure of SEM is modeled as a GGM:
\begin{equation}
\label{psychometrika:eq:RNM}
\hat{\pmb{\Sigma}} = \pmb{\Lambda}  \left( \pmb{I} - \pmb{B}   \right)^{-1} \pmb{\Psi} \left( \pmb{I} - \pmb{B}  \right)^{-1\top} \pmb{\Lambda}^{\top} + \pmb{\Delta}_{\pmb{\Theta}} \left( \pmb{I} - \pmb{\Omega}_{\pmb{\Theta}} \right)^{-1} \pmb{\Delta}_{\pmb{\Theta}}.
\end{equation}
This modeling framework conceptualizes latent variable and network modeling as two sides of the same coin, and offers immediate benefits to both. In latent variable modeling, RNM allows for the estimation of a factor structure (possibly including structural relationships between the latent variables), while having no uncorrelated errors and thus no local independence. The error-correlations, however, are still highly structured due to the residual network structure. This can be seen as a compromise between the ideas of network analysis and factor modeling; while we agree that local independence is plausibly violated in many psychometric tests, we think the assumption of no underlying latent traits and therefore a sparse GGM may often be too strict. For network modeling, RNM allows a researcher to estimate a sparse network structure while taking into account that some of the covariation between items was caused by a set of latent variables. Not taking this into account would lead to a saturated model \citep{chandrasekaran2010latent}, whereas the residual network structure can be sparse.

To avoid confusion between residual correlations, we will denote edges in $\pmb{\Omega}_{\pmb{\Theta}}$ \emph{residual interactions}. Residual interactions can be understood as pairwise linear effects, possibly due to some causal influence or partial overlap between indicators that is left after controlling for the latent structure. Consider again the indicators for agreeableness ``Am indifferent to the feelings of others'' and ``Inquire about others' well-being''. It seems clear that we would not expect these indicators to be locally independent after conditioning on agreeableness; being indifferent to the feelings of others will cause one to not inquire about other's well-being. Thus, we could expect these indicators to feature a residual interaction; some degree of correlation between these indicators is expected to remain, even after conditioning on the latent variable and all other indicators in the model. 

The RNM framework in particular offers a new way of improving the fit of confirmatory factor models. In contrast to increasingly popular methods such as exploratory SEM (ESEM; \citealt{marsh2014exploratory}) or LASSO regularized SEM models \citep{regsem}, the RNM framework improves the fit by adding residual interactions rather than allowing for more cross-loadings. The factor structure is kept exactly intact as specified in the confirmatory model. Importantly, therefore, the interpretation of the latent factor does not change. This can be highly valuable in the presence of a strong theory on the latent variables structure underlying a dataset even in the presence of violations of local independence.

\section{Exploratory Network Estimation}

Both the LNM and RNM modeling frameworks allow for confirmative testing of network structures. Confirmatory estimation is straightforward and similar to estimating SEM models, with the exception that instead of modeling $\pmb{\Psi}$ or $\pmb{\Theta}$ now the latent network $\pmb{\Omega}_{\pmb{\Psi}}$ or $\pmb{\Omega}_{\pmb{\Theta}}$ is modeled. Furthermore, both modeling frameworks allow for the confirmatory fit of a network model. In LNM, a confirmatory network structure can be tested by setting $\pmb{\Lambda}=\pmb{I}$ and $\pmb{\Theta} = \pmb{O}$; in RNM, a confirmatory network model can be tested by omitting any latent variables. We have developed the R package \textVerb{lvnet}\footnote{Developmental version: \url{http://www.github.com/sachaepskamp/lvnet}; stable version: \\ \url{https://cran.r-project.org/package=lvnet}}, which utilizes \textVerb{OpenMx} \citep{neale2015openmx} for confirmative testing of RNM and LNM models (as well as a combination of the two). The \textVerb{lvnet} function can be used for this purpose by specifying the fixed and the free elements of model matrices. The package returns model fit indices (e.g., the RMSEA, CFI and $\chi^2$ value), parameter estimates, and allows for model comparison tests. 

Often the network structure, either at the residual or the latent level, is unknown and needs to be estimated. To this end, the package includes two exploratory search algorithms described below: step-wise model search and penalized maximum likelihood estimation. For both model frameworks and both search algorithms, we present simulation studies to investigate the performance of these procedures. As is typical in simulation studies investigating the performance of network estimation techniques, we investigated the \emph{sensitivity} and \emph{specificity} \citep{van2014new}. These measures investigate the estimated edges versus the edges in the true model, with a `positive' indicating an estimated edge and a `negative' indicating an edge that is estimated to be zero. Sensitivity, also termed the true positive rate, gives the ratio of the number of true edges that were detected in the estimation versus the total number of edges in the true model:
\[
\text{sensitivity} = \frac{\text{\# true positives}}{\text{\# true positives} + \text{\# of false negatives}} 
\]
Specificity, also termed the true negative rate, gives the ratio of true missing edges detected in the estimation versus the total number of absent edges in the true model:
\[
\text{specificity} = \frac{\text{\# true negatives}}{\text{\# true negatives} + \text{\# false positives}} 
\]
The specificity can be seen as a function of the number of false positives: a high specificity indicates that there were not many edges detected to be nonzero that are zero in the true model. To favor degrees of freedom, model sparsity and interpretability, specificity should be high all-around---estimation techniques should not result in many false positives---whereas sensitivity should increase as a function of the sample size.

\subsection{Simulating Gaussian Graphical models}

In all simulation studies reported here, networks were constructed in the same way as done by \citet{yin2011sparse} in order to obtain a positive definite inverse-covariance matrix $\pmb{K}$. First, a network structure was generated without weights. Next, weights were drawn randomly from a uniform distribution between $0.5$ and $1$, and made negative with $50\%$ probability. The diagonal elements of $\pmb{K}$ were then set to $1.5$ times the sum of all absolute values in the corresponding row, or $1$ if this sum was zero. Next, all values in each row were divided by the diagonal value, ensuring that the diagonal values become $1$. Finally, the matrix was made symmetric by averaging the lower and upper triangular elements. In the chain graphs used in the following simulations, this algorithm created networks in which the non-zero partial correlations had a mean of $0.33$ and a standard deviation of $0.04$. 

\subsection{Stepwise Model Search}

In exploratory search, we are interested in recovering the network structure of either $\pmb{\Omega}_{\pmb{\Psi}}$ in LNM or $\pmb{\Omega}_{\pmb{\Theta}}$ in RNM. This can be done through a step-wise model search, either based on $\chi^2$ difference tests (Algorithm \ref{psychometrika:stepwiseChi}) or on minimization of some information criterion (Algorithm \ref{psychometrika:stepwiseIC}) such as the Akaike information criterion (AIC), Bayesian information criterion (BIC) or the extended Bayesian information criterion (EBIC; \citealt{chen2008EBIC}) which is now often used in network estimation \citep{van2014new, foygel2010extended}. In LNM, removing edges from $\pmb{\Omega}_{\pmb{\Psi}}$ cannot improve the fit beyond that of an already fitting CFA model. Hence, model search for $\pmb{\Omega}_{\pmb{\Psi}}$ should start at a fully populated initial setup for $\pmb{\Omega}_{\pmb{\Psi}}$. In RNM, on the other hand, a densely populated $\pmb{\Omega}_{\pmb{\Theta}}$ would lead to an over-identified model, and hence the step-wise model search should start at an empty network $\pmb{\Omega}_{\pmb{\Theta}} = \pmb{O}$. The function \textVerb{lvnetSearch} in the \textVerb{lvnet} package can be used for both search algorithms.

\begin{algorithm}[H]
\begin{algorithmic}
\STATE Start with initial setup for $\pmb{\Omega}$
\REPEAT
\FORALL{Unique elements of $\pmb{\Omega}$ }
\STATE{Remove edge if present or add edge if absent}
\STATE{Fit model with changed edge}
\ENDFOR
\IF{Adding an edge significantly improves fit ($\alpha = 0.05$)}
\STATE{Add edge that improves fit the most}
\ELSIF{Removing an edge does not significantly worsen fit ($\alpha = 0.05$)}
\STATE{Remove edge that worsens fit the least}
\ENDIF
\UNTIL{No added edge significantly improves fit and removing any edge significantly worsens fit}
\end{algorithmic}
 \caption{Stepwise network estimation by $\chi^2$ difference testing.}
 \label{psychometrika:stepwiseChi}
\end{algorithm}

\begin{algorithm}[H]
\begin{algorithmic}
\STATE Start with initial setup for $\pmb{\Omega}$
\REPEAT
\FORALL{Unique elements of $\pmb{\Omega}$ }
\STATE{Remove edge if present or add edge if absent}
\STATE{Fit model with changed edge}
\ENDFOR
\IF{Any changed edge improved AIC, BIC or EBIC}
\STATE{Change edge that improved AIC, BIC or EBIC the most}
\ENDIF
\UNTIL{No changed edge improves AIC, BIC or EBIC}
\end{algorithmic}
 \caption{Stepwise network estimation by AIC, BIC or EBIC optimization.}
 \label{psychometrika:stepwiseIC}
\end{algorithm}

\subsubsection{Simulation Study 1: Latent Network Modeling}

\begin{figure}
\centering
\includegraphics[width = 1\textwidth]{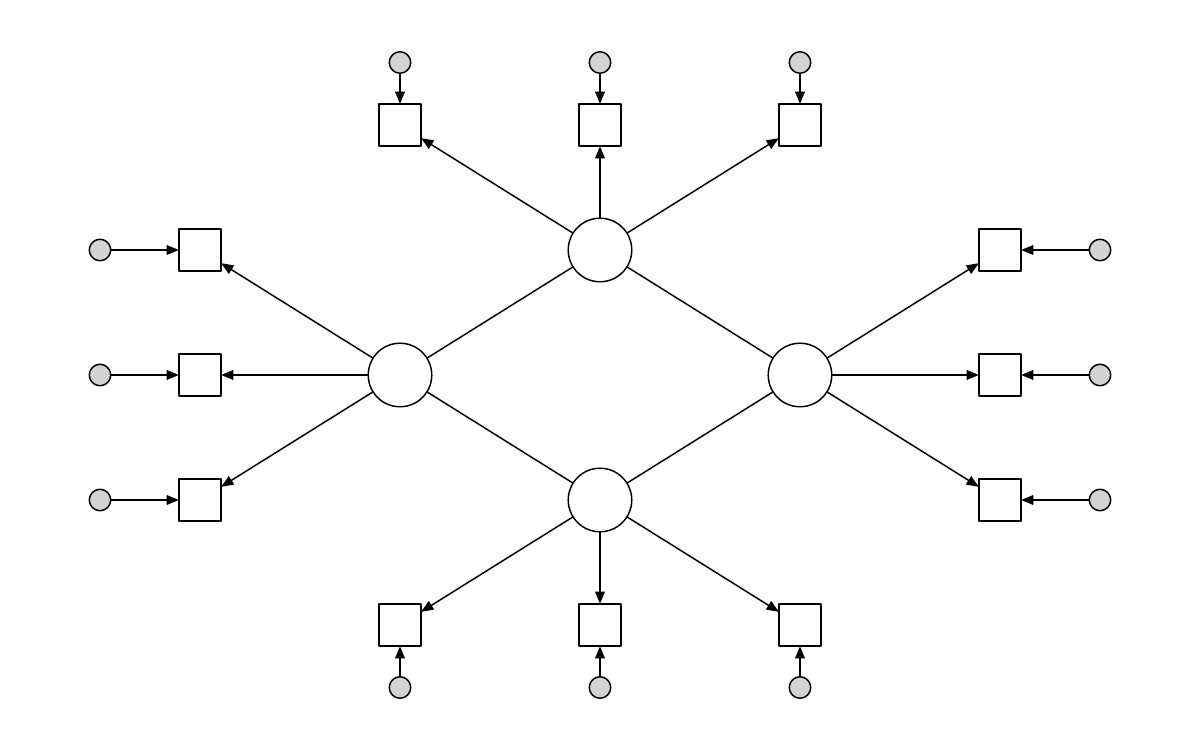}
\caption{Model used in simulation study 1: step-wise model search in latent network modeling. Four latent variables were each measured by three items. Latent variables covary due to the structure of a latent Gaussian graphical model in which edges indicate partial correlation coefficients. This model has the form of a chain graph, which cannot be represented in a structural equation model. Factor loadings, residual variances and latent variances were set to $1$ and the latent partial correlations had an average of $0.33$ with a standard deviation of $0.04$.}
\label{psychometrika:fig:LNMsim}
\end{figure}

We performed a simulation study to assess the performance of the above mentioned step-wise search algorithms in LNM models. Figure~\ref{psychometrika:fig:LNMsim} shows the LNM model under which we simulated data. In this model, four latent factors with three indicators each were connected in a latent network. The latent network was a chain network, leading all latent variables to be correlated according to a structure that cannot be represented in SEM. Factor loadings and residual variances were set to 1, and the network weights were simulated as described in the section ``Simulating Gaussian Graphical models''. The simulation study followed a $5 \times 4$ design: the sample size was varied between 50, 100, 250, 500 and $1\,000$ to represent typical sample sizes in psychological research, and the stepwise evaluation criterion was either $\chi^2$ difference testing, AIC, BIC or EBIC (using a tuning parameter of 0.5). Each condition was simulated $1\,000$ times, resulting in $20\,000$ total simulated datasets.

\begin{figure}
    \centering
        \includegraphics[width=\textwidth,page=1]{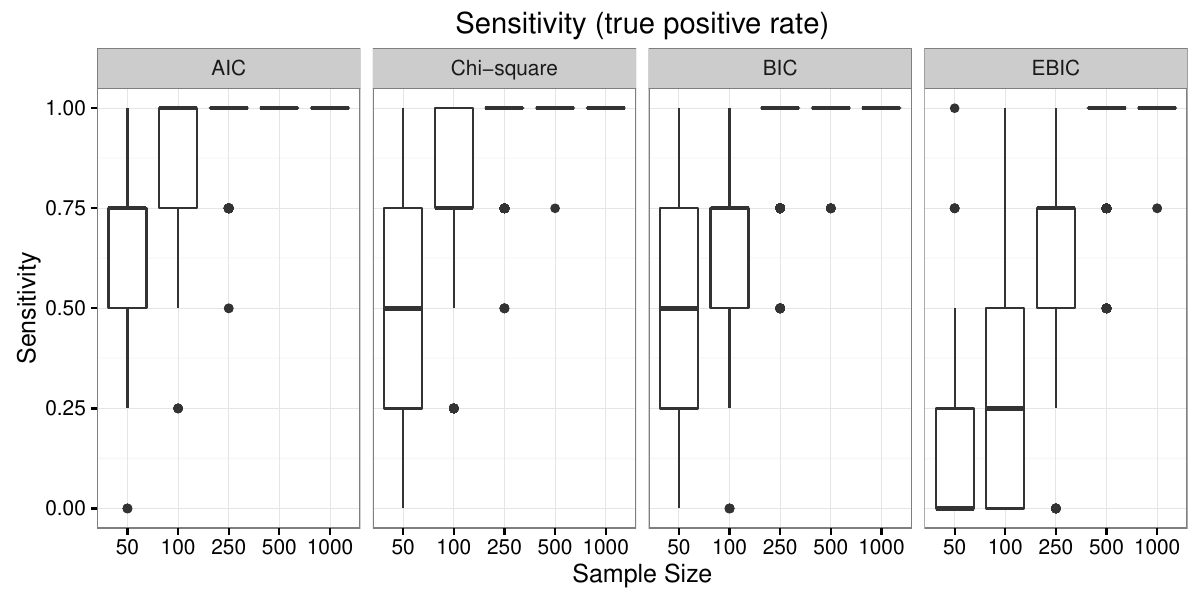}
        \includegraphics[width=\textwidth,page=2]{Figure5.pdf}
    \caption{Simulation results of simulation study 1: step-wise model search in latent network modeling. Each condition was replicated $1\,000$ times, leading to $20\,000$ total simulated datasets. High sensitivity indicates that the method is able to detect edges in the true model, and high specificity indicates that the method does not detect edges that are zero in the true model.}
    \label{psychometrika:fig:LNMresults}
\end{figure}

Figure~\ref{psychometrika:fig:LNMresults} shows the results of the simulation study. Data is represented in standard boxplots \citep{mcgill1978variations}: the box shows the 25th, 50th (median) 75th quantiles, the whiskers range from the largest values in $1.5$ times the inter-quantile range (75th - 25th quantile) and points indicate outliers outside that range. In each condition, we investigated the sensitivity and specificity. The top panel shows that sensitivity improves with sample size, with AIC performing best and EBIC worst. From sample sizes of 500 and higher all estimation criterion performed well in retrieving the edges. The bottom panel shows that specificity is generally very high, with EBIC performing best and AIC worst. These results indicate that the step-wise procedure is conservative and prefers simpler models to more complex models; missing edges are adequately detected but present edges in the true model might go unnoticed except in larger samples. With sample sizes over 500, all four estimation methods show both a high sensitivity and specificity.

\subsubsection{Simulation Study 2: Residual Network Modeling}

\begin{figure}
    \centering
        \includegraphics[width=\textwidth,page=1]{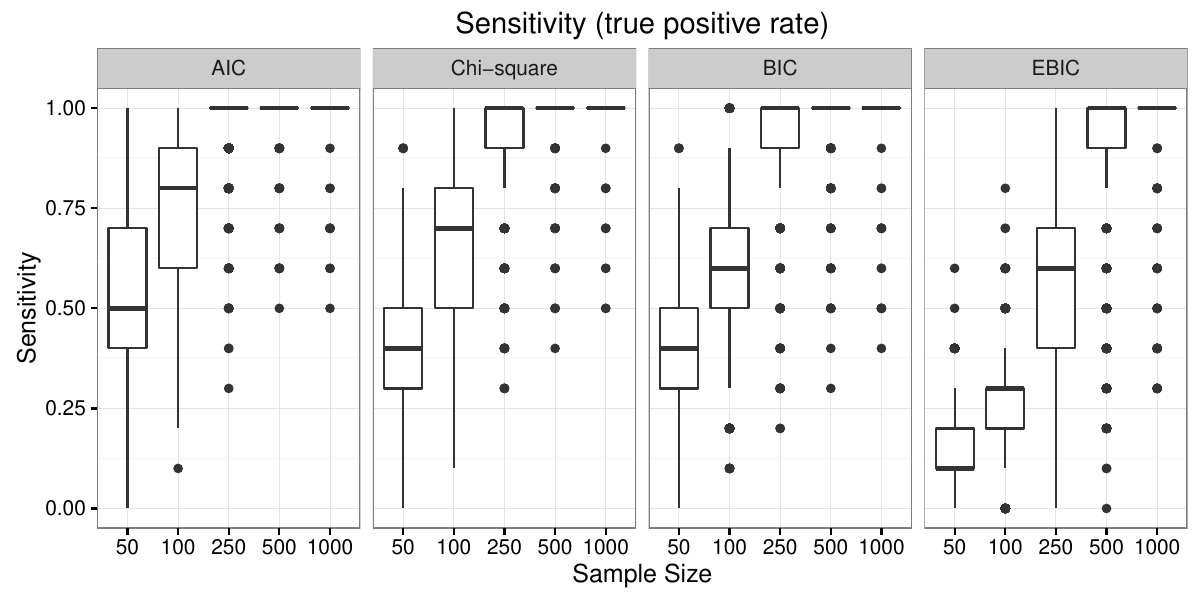}
        \includegraphics[width=\textwidth,page=2]{Figure7.pdf}
    \caption{Simulation results of simulation study 2: step-wise model search in residual network modeling. Each condition was replicated $1\,000$ times, leading to $20\,000$ total simulated datasets. High sensitivity indicates that the method is able to detect edges in the true model, and high specificity indicates that the method does not detect edges that are zero in the true model.}
    \label{psychometrika:fig:RNMresults}
\end{figure}

We conducted a second simulation study to assess the performance of step-wise model selection in RNM models. Figure~\ref{psychometrika:fig:RNMresults} shows the model under which data were simulated: two latent variables with 5 indicators each. The residual network was constructed to be a chain graph linking a residual of an indicator of one latent variable to two indicators of the other latent variable. This structure cannot be represented by a DAG and causes all residuals to be connected, so that $\pmb{\Theta}$ is fully populated. Factor loadings and residual variances were set to 1, the factor covariance was set to $0.25$, and the network weights were simulated as described in the section ``Simulating Gaussian Graphical models''.

\begin{figure}
    \centering
        \includegraphics[width = 1\textwidth]{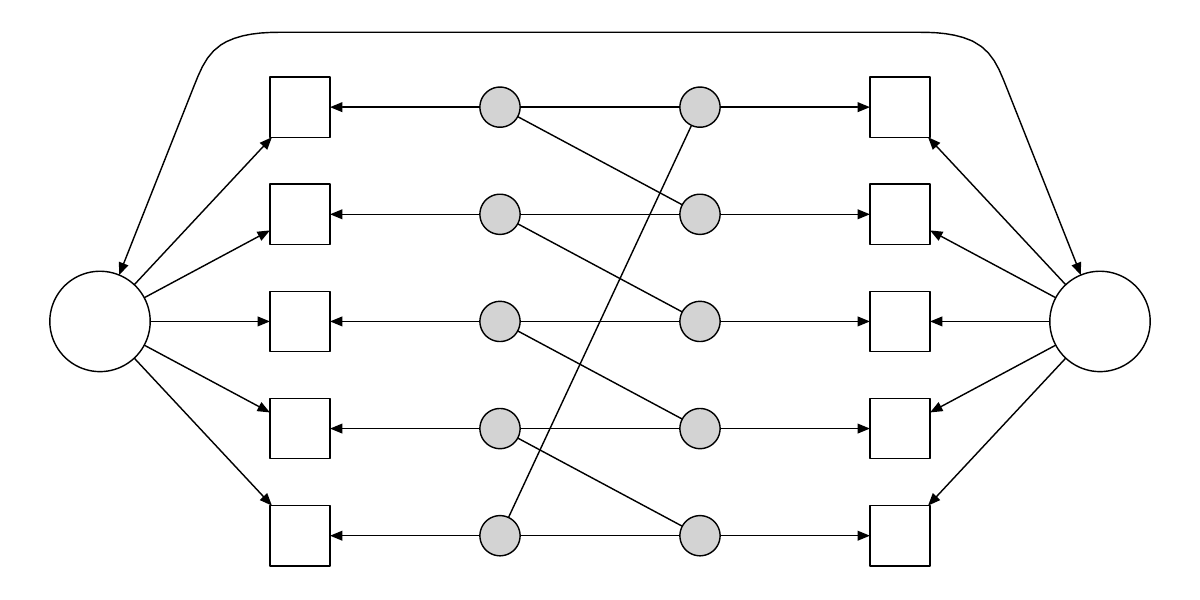}
    \caption{Model used in simulation study 2: step-wise model search in residual network modeling. Two latent variables were each measured by five items; a Gaussian graphical model, in which edges indicate partial correlation coefficients, leads to all residuals to be correlated due to a chain graph between residuals, which cannot be represented in a structural equation model. Factor loadings, residual variances and latent variances were set to $1$, the factor covariance was set to $0.25$ and the latent partial correlations had an average of $0.33$ with a standard deviation of $0.04$.}
    \label{psychometrika:fig:RNMsimmodels}
\end{figure}

The simulation study followed a $5 \times 4$ design; sample size was again varied between 50, 100, 250, 500 and $1\,000$, and models were estimated using either $\chi^2$ significance testing, AIC, BIC or EBIC. Factor loadings and factor variances were set to $1$ and the factor correlation was set to $0.25$. The weights in $\pmb{\Omega}_{\pmb{\Theta}}$ were chosen as described in the section ``Simulating Gaussian Graphical models''. Each condition was simulated $1,000$ times, leading to $20\,000$ total datasets.

Figure~\ref{psychometrika:fig:RNMresults} shows the results of the simulation study. The top panel shows that sensitivity increases with sample size and performs best when using AIC as the criterion. BIC performed comparably in sensitivity to $\chi^2$ testing and EBIC performed the worst. The bottom panel shows that specificity was very high for all sample sizes and all criteria, with EBIC performing best and AIC worst. These results indicate that the number of false positives is very low and that the method is on average well capable of discovering true edges for sample size larger than 250. In sum, all four criteria perform well with EBIC erring on the side of caution and AIC erring on the side of discovery.

\subsection{LASSO Regularization}

While the step-wise model selection algorithms perform well in retrieving the correct network structure, they are very slow when the number of nodes in the network increases (e.g., more than $10$ nodes). This is particularly important in the context of RNM, in which the number of indicators can be larger than 10 even in small models. A popular method for fast estimation of high-dimensional network structures is by applying the least absolute shrinkage and selection operator (LASSO; \citealt{tibshirani1996regression}). LASSO regularization has also recently been introduced in the SEM literature \citep{regsem} as a method for obtaining sparser structures of $\pmb{\Lambda}$ and $\pmb{B}$. In the LASSO, instead of optimizing the likelihood function as described in \eqref{psychometrika:loglik}, a penalized likelihood is optimized \citep{regsem}:
\begin{equation}
\label{psychometrika:penloglik}
\min_{\hat{\pmb{\Sigma}}} \left[ \log \det \left( \hat{\pmb{\Sigma}}\right)   + \mathrm{Trace}\left( \pmb{S} \hat{\pmb{\Sigma}}^{-1} \right)  - \log \det \left( \hat{\pmb{S}}\right)  - P  + \nu \mathrm{Penalty} \right],
\end{equation}
in which $\nu$ denotes a tuning parameter controlling the level of penalization. The penalty here is taken to be the sum of absolute parameters:
\[
\mathrm{Penalty} = \sum_{<i, j>} |\omega_{ij}|,
\]
in which $\omega_{ij}$ denotes an element from either $\pmb{\Omega}_{\pmb{\Psi}}$ or $\pmb{\Omega}_{\pmb{\Theta}}$. Other penalty functions may be used as well---such as summing the squares of parameter estimates (ridge regression; \citealt{hoerl1970ridge}) or combining both absolute and squared values (elastic net; \citealt{zou2005regularization})---but these are not currently implemented in \textVerb{lvnet}. The benefit of the LASSO is that it returns models that perform better in cross-validation. In addition, the LASSO yields sparse models in which many relationships are estimated to be zero.

\begin{algorithm}[H]
\begin{algorithmic}
\FORALL{Sequence of tuning parameters $\nu_1, \nu_2, \ldots$}
\STATE{Estimate LASSO regularized model using given tuning parameter}
\STATE{Count the number of parameters for which the absolute estimate is larger than $\epsilon$}
\STATE{Determine information criterion (AIC, BIC or EBIC) given fit and number of parameters}
\ENDFOR
\STATE{Select model with best AIC, BIC or EBIC}
\STATE{Refit this model without LASSO in which absolute parameters smaller than $\epsilon$ are fixed to zero}
\end{algorithmic}
 \caption{LASSO estimation for exploratory network search.}
 \label{psychometrika:LASSOalg}
\end{algorithm}

The \textVerb{lvnet} function allows for LASSO regularization for a given model matrix ($\pmb{\Omega}_{\pmb{\Theta}}$, $\pmb{\Omega}_{\pmb{\Psi}}$, $\pmb{\Theta}$, $\pmb{\Psi}$, $\pmb{\Lambda}$ or $\pmb{B}$) and a given value for the tuning parameter $\nu$. The optimizer used in \textVerb{lvnet} does not return exact zeroes. To circumvent this issue, any absolute parameter below some small value $\epsilon$ (by default $\epsilon = 0.0001$) is treated as zero in counting the number of parameters and degrees of freedom \citep{zou2007degrees}. The \textVerb{lvnetLasso} function implements the search algorithm described in Algorithm \ref{psychometrika:LASSOalg} to automatically choose an appropriate tuning parameter, use that for model selection and rerun the model to obtain a comparable fit to non-regularized models. In this algorithm, a sequence of tuning parameters is tested, which is set by default to a logarithmically spaced sequence of 20 values between $0.01$ and $1$.

\subsubsection{Simulation Study 3: Latent Network Modeling}

\begin{figure}
    \centering
        \includegraphics[width=\textwidth,page=1]{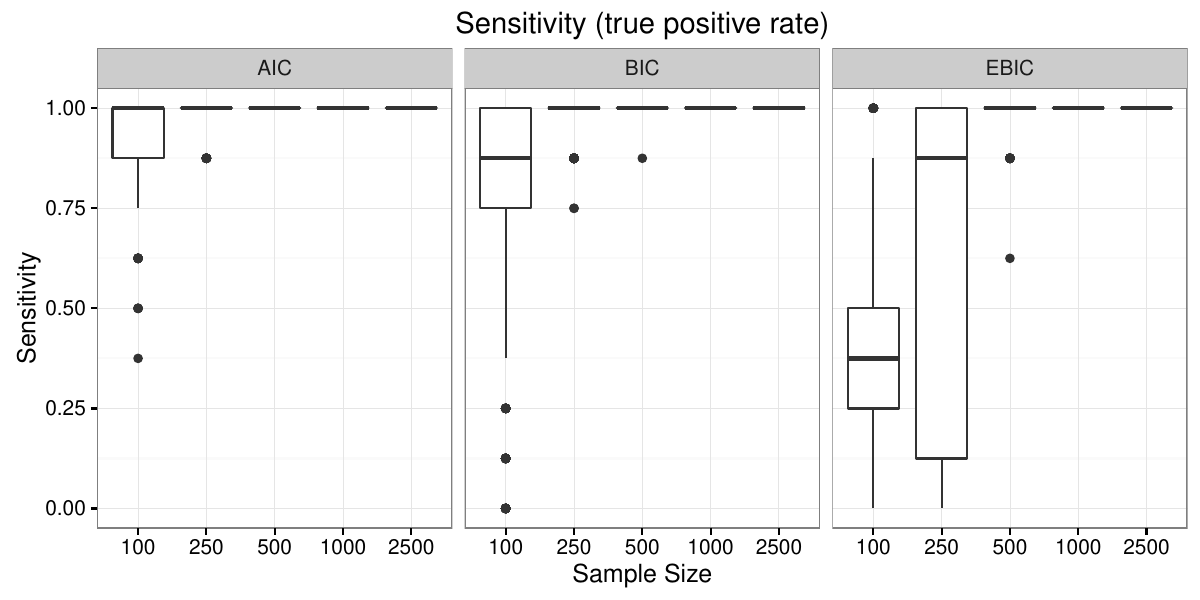}
        \includegraphics[width=\textwidth,page=2]{Figure8.pdf}
    \caption{Simulation results of simulation study 3: model selection via penalized maximum likelihood estimation in latent network modeling. The same model as in Figure~\ref{psychometrika:fig:LNMsim} was used except now with $4$ latent variables leading to $24$ observed variables. For each sample size $1\,000$ datasets were generated, leading to $5\,000$ total simulated datasets on which AIC, BIC or EBIC was used to select the best model. High sensitivity indicates that the method is able to detect edges in the true model, and high specificity indicates that the method does not detect edges that are zero in the true model.}
    \label{psychometrika:fig:LNM_lassoresults}
\end{figure}

We studied the performance of LASSO penalization in estimating the latent network structure in a similar simulation study to the study of the step-wise procedure described above. Data were simulated under a similar model to the one shown in Figure~\ref{psychometrika:fig:LNMsim}, except that now $8$ latent variables were used leading to a total of $24$ observed variables. All parameter values were the same as in simulation study 1. The simulation followed a $5 \times 3$ design. Sample size was varied between $100$, $250$, $500$, $1\,000$ and $2\,500$, and for each sample size $1\,000$ datasets were simulated leading to a total of $5\,000$ generated datasets. On these datasets the best model was selected using either AIC, BIC or EBIC, leading to $15\,000$ total replications. In each replication, sensitivity and specificity were computed. Figure~\ref{psychometrika:fig:LNM_lassoresults} shows that AIC had a relatively poor specificity all-around, but a high sensitivity. EBIC performed well with sample sizes of $500$ and higher.

\subsubsection{Simulation Study 4: Residual Network Modeling}

\begin{figure}
    \centering
        \includegraphics[width=\textwidth,page=1]{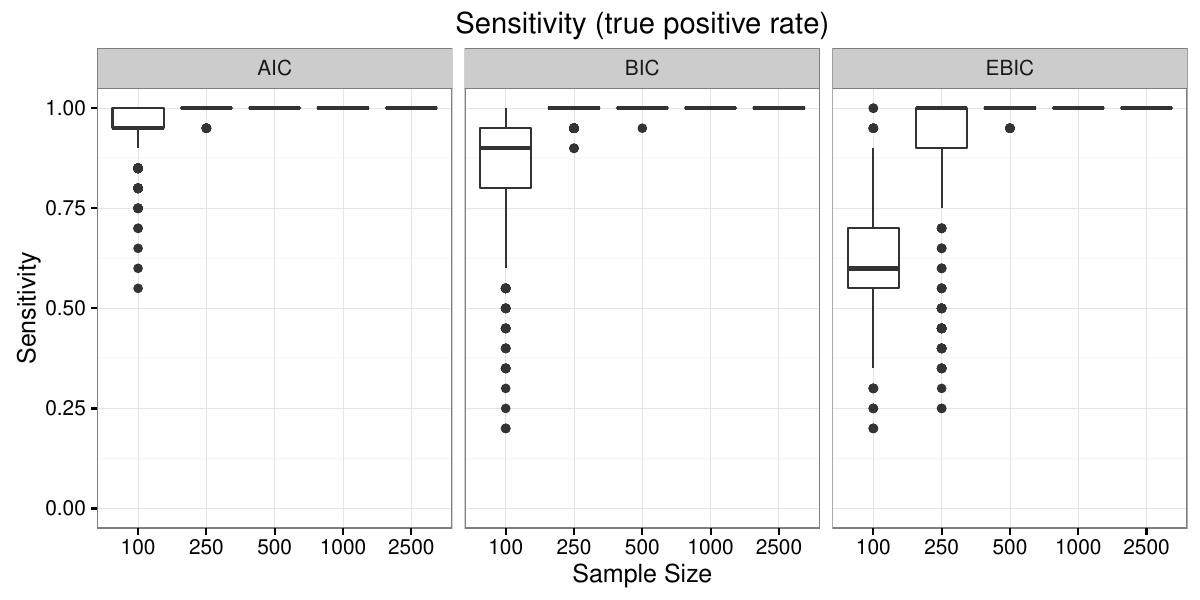}
        \includegraphics[width=\textwidth,page=2]{Figure9.pdf}
    \caption{Simulation results of simulation study 4: model selection via penalized maximum likelihood estimation in residual network modeling. The same model as in Figure~\ref{psychometrika:fig:RNMsimmodels} was used except now with $4$ latent variables leading to $20$ observed variables. For each sample size $1\,000$ datasets were generated, leading to $5\,000$ total simulated datasets on which AIC, BIC or EBIC was used to select the best model. High sensitivity indicates that the method is able to detect edges in the true model, and high specificity indicates that the method does not detect edges that are zero in the true model.}
    \label{psychometrika:fig:RNM_lassoresults}
\end{figure}

To assess the performance of LASSO in estimating the residual network structure we simulated data as in Figure~\ref{psychometrika:fig:RNMsimmodels}, except that in this case four latent variables were used, each with 5 indicators, the residuals of which were linked via a chain graph. All parameter values were the same as in simulation study 2.  The design was the same as in simulation study 3, leading to $5\,000$ generated datasets on which AIC, BIC or EBIC were used to select the best model. While Figure~\ref{psychometrika:fig:RNM_lassoresults} shows good performance of the LASSO in retrieving the residual network structure and similar results as before: AIC performs the worst in specificity and EBIC the best.

\section{Empirical Example: Personality Inventory}

In this section, we demonstrate LNM and RNM models by confirmatively testing a model and exploratively searching a residual and a latent network structure. We use the \textVerb{lvnet} package, which can be installed in R as follows:
\begin{verbatim}
> install.packages("lvnet")
> library("lvnet")
\end{verbatim}

To exemplify the method, we will use a dataset from the psych package \citep{psych} on the Big 5 personality traits \citep{benet1998cinco, digman1989five,goldberg2001alternative,goldberg1998structure, mccrae1997personality}. This dataset consists of 2800 observations of 25 items designed to measure the 5 central personality traits with 5 items per trait. We estimated the CFA model on the BFI dataset. Next, we used LASSO estimation to the RNM model using $100$ different tuning parameters and using EBIC as criterion to maximize specificity and search for a sparse residual network. The fully correlated covariance matrix of latent variables is equivalent to a fully connected latent network structure. Thus, after fitting a RNM model, we can again apply LASSO to estimate a latent network in the resulting model, which we abbreviate here to an RNM+LNM model. The R code used for this analysis can be found in the supplementary materials.

\begin{table}[ht]
\centering
{\footnotesize
\begin{tabular}{rrrrrrrrr}
  \hline
 & df & $\chi^2$ & AIC & BIC & EBIC & RMSEA & TLI & CFI \\ 
  \hline
CFA & 265 & 4713.94 & 183233.7 & 183589.9 & 184542.4 & 0.08 & 0.75 & 0.78 \\ 
  RNM & 172 & 806.63 & 179511.0 & 180419.4 & 182848.2 & 0.04 & 0.94 & 0.97 \\ 
  RNM+LNM & 176 & 843.18 & 179539.5 & 180424.2 & 182789.5 & 0.04 & 0.94 & 0.97 \\
   \hline
\end{tabular}
}
\caption{Fit measures for three models estimated on the BFI dataset in the psych R package. CFA is the correlated five-factor model. RNM is the same model as the CFA model with a residual network. RNM+LNM denotes the same model as the RNM model in which edges of the latent network have been removed.}
\label{psychometrika:tab:bfi}
\end{table}

\begin{figure}
    \centering
        \includegraphics[width=\textwidth,page=1]{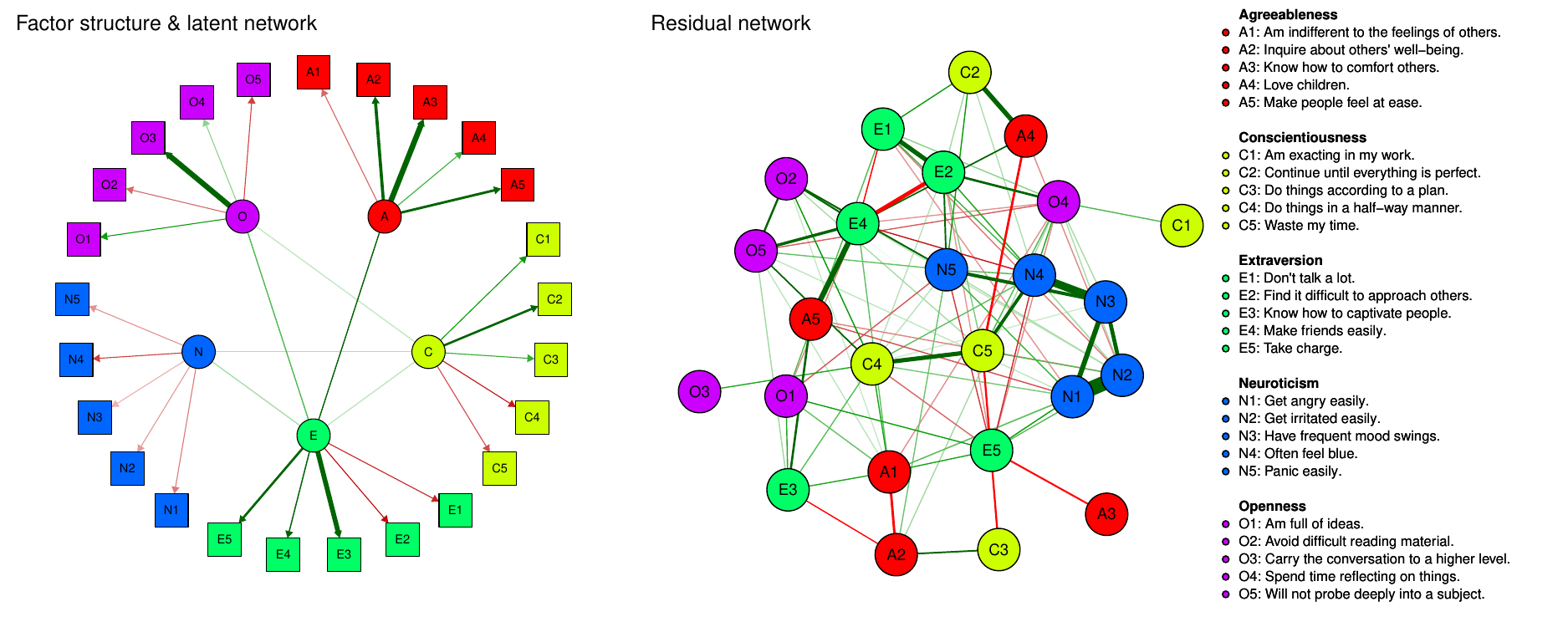}
    \caption{Visualization of the factor structure and latent network (left) and the residual network (right) of the BFI personality dataset from the  psych package in R. LASSO estimation with $100$ different tuning parameters in combination with EBIC model selection was used to first estimate the residual network structure and following the latent network structure.}
    \label{psychometrika:fig:BFI}
\end{figure}

Table~\ref{psychometrika:tab:bfi} shows the fit of the three models. The CFA model fits poorly. The RNM model has substantively improved fit and resulted in good fit indices overall. The estimated RNM+LNM model showed that $5$ edges could be removed from the latent network after taking residual interactions into account. Figure~\ref{psychometrika:fig:BFI} shows the factor structure and residual network of the final RNM+LNM model. It can be seen that Agreeableness is now only connected to extraversion: after taking into account someone's level of extraversion agreeableness is independent of the other three personality traits. Extraversion is the most central node in this network and the only trait that is directly linked to all other traits. The residual network shows many meaningful connections. While seemingly densely connected, this network only has 30\% of all possible edges in a network of that size, leading the model to have $176$ degrees of freedom. The corresponding residual covariance structure is fully populated with no zero elements.

It should be noted that the procedures used in this example are highly explorative. While Table~\ref{psychometrika:tab:bfi} shows that the RNM fits better than the CFA model, the CFA model is solely based on theory whereas the RNM model was found through high-dimensional model search. As a result, to substantially interpret the structure found in Figure~\ref{psychometrika:fig:BFI} it should first be replicated in independent samples.

\section{Conclusion}

In this paper we introduced a formal psychometric model for network modeling of multivariate normal data. We contrasted this model with latent variable models as commonly used in CFA and SEM. Furthermore, using the CFA and SEM frameworks, we proposed two generalizations of the network model to encompass latent variable structures within the network paradigm. In the first generalization, LNM, we construct a network among the latent variables, whereas in the second generalization, RNM, a network is formed among residuals of indicators. Both frameworks offer powerful benefits over both latent variable and network modeling. From the perspective of latent variable modeling, the LNM framework allows one to exploratively search for conditional independence relationships between latent variables without the need for prior theory, and the RNM framework allows one to model latent common causes without assuming local independence. From the perspective of network modeling, the LNM framework allows one to model network structures while taking measurement error into account, and the RNM framework allows one to estimate a network structure, even when all nodes are in part caused by unobserved or latent variables. 
In addition, both frameworks allow for network models to be fitted and compared to SEM models. The discussed methodology has been implemented in the freely available R package \textVerb{lvnet}.

Simulation studies showed that step-wise search and penalized maximum likelihood estimation of the residual or latent network structures resulted in high specificity all around---the methods did not result often in false positives---and rapidly increasing sensitivity as a function of the sample size; the higher the sample size, the more true edges were detected in the algorithm. These numbers are comparable to state-of-the-art network estimation techniques in sample and model sizes that are plausible in psychological settings \citep{van2014new}. In all four simulation studies, using AIC as the model selection criterion led to the best sensitivity and using EBIC led to the best specificity. However, it is important to note that the choice of a particular information criterion cannot be argued by these numbers alone, and depends on the relative importance one assigns to the side of discovery (optimizing sensitivity) or the side of caution (optimizing specificity; \citealt{dziak2012sensitivity}). Furthermore, it should be noted that these simulation results are specific for the particular model setup and sample sizes used; results might be different for other kinds of models or sample size ranges.

In addition to the LNM and RNM frameworks, other combinations of CFA, SEM and network modeling are possible as well. For example, a framework can be constructed which contains both a latent and a residual network (as shown in our empirical example), or directed regression paths as in the SEM model can be added to the LNM model. While these models are all estimable in the \textVerb{lvnet} software, in the current paper we chose to focus on the distinct benefits that modeling a residual or latent network presents. Thus, in this paper, we only described the modeling of multivariate normal data. More advanced models are possible, but not yet implemented in the \textVerb{lvnet} software. In the case of binary variables, the appropriate model to use is the Ising model, which has been shown to be equivalent to multivariate item response models \citep{netpsych,marsman2015bayesian}. Future research could aim at constructing Ising models among binary latent variables in latent class analysis, or constructing residual networks in models with binary indicators. Finally, the expressions optimized in Equations \eqref{psychometrika:loglik} and \eqref{psychometrika:penloglik} are based on summary statistics and therefore only truly applicable to complete data. With incomplete data, the appropriate estimation method is to use full-information maximum likelihood (FIML; \citealt{arbuckle1996full}); however, FIML has not yet been implemented in the \textVerb{lvnet} software.  

In our view, the presented modeling framework is a versatile and promising addition to the spectrum of psychometric models. The GGM, which has a central place in this modeling framework, acts as a natural interface between correlation and causality, and we think this representation should receive more attention in psychometrics. From the point of view afforded by the current paper, the typical attempt to determine directed SEMs from correlation structures in fact appears somewhat haphazard in psychology, a historical accident in a field that has been prematurely directed to hypothesis testing at the expense of systematic exploration. Perhaps, psychometrics as a field should consider taking a step back to focus on the consistent identification of GGMs, instead of wanting to jump to the causal conclusion immediately. In this regard, the fact that GGMs do not have equivalent models would appear to be a major benefit, as they allow us to focus on charting connections between variables systematically, without being forced to adhere to one particular causal interpretation or another. In addition, because the GGM does not specify the nature or direction of interactions between variables, it appears a natural model for research situations where no temporal information or experimental interventions are present, so that associations may arise for a multitude of reasons: the GGM can be consistently interpreted regardless of whether associations arise from direct causal relations, reciprocal causation, latent common causes, semantic overlap between items, or homeostatic couplings of parameters. While this can be seen as a downside of the GGM---the lack of directionality leads to fewer falsifiable hypotheses---it can also be a major asset in a field like psychology, where strong causal theory is sparse and the identification of DAGs often appears a bridge too far.

\bibliographystyle{apalike}
\bibliography{Bibliography}

\begin{thebibliography}{}

\bibitem[Arbuckle et~al., 1996]{arbuckle1996full}
Arbuckle, J.~L., Marcoulides, G.~A., and Schumacker, R.~E. (1996).
\newblock Full information estimation in the presence of incomplete data.
\newblock {\em Advanced structural equation modeling: Issues and techniques},
  pages 243--277.

\bibitem[Benet-Martinez and John, 1998]{benet1998cinco}
Benet-Martinez, V. and John, O. (1998).
\newblock {Los Cinco Grandes across Cultures and Ethnic Groups: Multitrait
  Multimethod Analyses of the Big Five in Spanish and English}.
\newblock {\em Journal of Personality and Social Psychology}, 75:729--750.

\bibitem[Bentler, 1990]{bentler1990comparative}
Bentler, P.~M. (1990).
\newblock Comparative fit indexes in structural models.
\newblock {\em Psychological bulletin}, 107(2):238--346.

\bibitem[Borsboom, 2008]{borsboom2008psychometric}
Borsboom, D. (2008).
\newblock {Psychometric Perspectives on Diagnostic Systems}.
\newblock {\em Journal of clinical psychology}, 64(9):1089--1108.

\bibitem[Borsboom and Cramer, 2013]{borsboom2013network}
Borsboom, D. and Cramer, A. O.~J. (2013).
\newblock Network analysis: An integrative approach to the structure of
  psychopathology.
\newblock {\em Annual Review of Clinical Psychology}, 9:91--121.

\bibitem[Borsboom et~al., 2011]{borsboom2011small}
Borsboom, D., Cramer, A. O.~J., Schmittmann, V.~D., Epskamp, S., and Waldorp,
  L.~J. (2011).
\newblock The small world of psychopathology.
\newblock {\em PloS one}, 6(11):e27407.

\bibitem[Browne and Cudeck, 1992]{browne1992alternative}
Browne, M.~W. and Cudeck, R. (1992).
\newblock Alternative ways of assessing model fit.
\newblock {\em Sociological Methods \& Research}, 21(2):230--258.

\bibitem[Chandrasekaran et~al., 2012]{chandrasekaran2010latent}
Chandrasekaran, V., Parrilo, P.~A., and Willsky, A.~S. (2012).
\newblock Latent variable graphical model selection via convex optimization
  (with discussion).
\newblock {\em The Annals of Statistics}, 40(4):1935--1967.

\bibitem[Chen and Chen, 2008]{chen2008EBIC}
Chen, J. and Chen, Z. (2008).
\newblock Extended bayesian information criteria for model selection with large
  model spaces.
\newblock {\em Biometrika}, 95(3):759--771.

\bibitem[Costantini et~al., 2015]{costantini2015state}
Costantini, G., Epskamp, S., Borsboom, D., Perugini, M., M\~{o}ttus, R.,
  Waldorp, L.~J., and Cramer, A. O.~J. (2015).
\newblock State of the {aRt} personality research: A tutorial on network
  analysis of personality data in {R}.
\newblock {\em Journal of Research in Personality}, 54:13--29.

\bibitem[Cramer et~al., 2012]{cramer2012dimensions}
Cramer, A. O.~J., Sluis, S., Noordhof, A., Wichers, M., Geschwind, N., Aggen,
  S.~H., Kendler, K.~S., and Borsboom, D. (2012).
\newblock Dimensions of normal personality as networks in search of
  equilibrium: You can't like parties if you don't like people.
\newblock {\em European Journal of Personality}, 26(4):414--431.

\bibitem[Cramer et~al., 2010]{cramer2010comorbidity}
Cramer, A. O.~J., Waldorp, L., van~der Maas, H., and Borsboom, D. (2010).
\newblock {Comorbidity: A Network Perspective}.
\newblock {\em Behavioral and Brain Sciences}, 33(2-3):137--150.

\bibitem[Dalege et~al., 2016]{dalege2016toward}
Dalege, J., Borsboom, D., van Harreveld, F., van~den Berg, H., Conner, M., and
  van~der Maas, H. L.~J. (2016).
\newblock {Toward a formalized account of attitudes: The Causal Attitude
  Network (CAN) model.}
\newblock {\em Psychological review}, 123(1):2--22.

\bibitem[Digman, 1989]{digman1989five}
Digman, J. (1989).
\newblock {Five Robust Trait Dimensions: Development, Stability, and Utility}.
\newblock {\em Journal of Personality}, 57(2):195--214.

\bibitem[Dziak et~al., 2012]{dziak2012sensitivity}
Dziak, J.~J., Coffman, D.~L., Lanza, S.~T., and Li, R. (2012).
\newblock Sensitivity and specificity of information criteria.
\newblock {\em The Methodology Center and Department of Statistics, Penn State,
  The Pennsylvania State University}.

\bibitem[Ellis and Junker, 1997]{ellis1997tail}
Ellis, J.~L. and Junker, B.~W. (1997).
\newblock Tail-measurability in monotone latent variable models.
\newblock {\em Psychometrika}, 62(4):495--523.

\bibitem[Epskamp et~al., ress]{netpsych}
Epskamp, S., Maris, G., Waldorp, L., and Borsboom, D. (in press).
\newblock Network psychometrics.
\newblock In Irwing, P., Hughes, D., and Booth, T., editors, {\em Handbook of
  Psychometrics}. Wiley, New York, NY, USA.

\bibitem[Epskamp et~al., 2016]{dynamics}
Epskamp, S., Waldorp, L.~J., M\~{o}ttus, R., and Borsboom, D. (2016).
\newblock Discovering psychological dynamics in time-series data.
\newblock {\em arXiv preprint}, page arXiv:1609.04156.

\bibitem[Foygel and Drton, 2010]{foygel2010extended}
Foygel, R. and Drton, M. (2010).
\newblock Extended {Bayesian} information criteria for {Gaussian} graphical
  models.
\newblock {\em Advances in Neural Information Processing Systems},
  23:2020--2028.

\bibitem[Fried et~al., 2015]{fried2015}
Fried, E.~I., Bockting, C., Arjadi, R., Borsboom, D., Amshoff, M., Cramer,
  O.~J., Epskamp, S., Tuerlinckx, F., Carr, D., and Stroebe, M. (2015).
\newblock {From loss to loneliness: The relationship between bereavement and
  depressive symptoms}.
\newblock {\em Journal of abnormal psychology}, 124(2):256--265.

\bibitem[Fried and van Borkulo, 2016]{Fried_van_Borkulo_2016}
Fried, E.~I. and van Borkulo, C. (2016).
\newblock Mental disorders as networks of problems: a review of recent
  insights.

\bibitem[Gates and Molenaar, 2012]{gates2012group}
Gates, K.~M. and Molenaar, P.~C. (2012).
\newblock Group search algorithm recovers effective connectivity maps for
  individuals in homogeneous and heterogeneous samples.
\newblock {\em NeuroImage}, 63(1):310--319.

\bibitem[Goldberg, 1993]{goldberg1998structure}
Goldberg, L. (1993).
\newblock {The Structure of Phenotypic Personality Traits}.
\newblock {\em American Psychologist}, 48(1):26--34.

\bibitem[Goldberg, 1990]{goldberg2001alternative}
Goldberg, L.~R. (1990).
\newblock An alternative ``description of personality'': the big-five factor
  structure.
\newblock {\em Journal of personality and social psychology}, 59(6):1216--1229.

\bibitem[Hayduk, 1987]{hayduk1987structural}
Hayduk, L.~A. (1987).
\newblock {\em Structural equation modeling with {LISREL}: Essentials and
  advances}.
\newblock Johns Hopkins University Press, Baltimore, MD, USA.

\bibitem[Hoerl and Kennard, 1970]{hoerl1970ridge}
Hoerl, A.~E. and Kennard, R.~W. (1970).
\newblock Ridge regression: Biased estimation for nonorthogonal problems.
\newblock {\em Technometrics}, 12(1):55--67.

\bibitem[Holland and Rosenbaum, 1986]{holland1986conditional}
Holland, P.~W. and Rosenbaum, P.~R. (1986).
\newblock Conditional association and unidimensionality in monotone latent
  variable models.
\newblock {\em The Annals of Statistics}, 14:1523--1543.

\bibitem[Ising, 1925]{ising1925beitrag}
Ising, E. (1925).
\newblock Beitrag zur theorie des ferromagnetismus.
\newblock {\em Zeitschrift f{\"u}r Physik A Hadrons and Nuclei},
  31(1):253--258.

\bibitem[Isvoranu et~al., 2016a]{isvoranua}
Isvoranu, A.~M., Borsboom, D., van Os, J., and Guloksuz, S. (2016a).
\newblock {A Network Approach to Environmental Impact in Psychotic Disorders:
  Brief Theoretical Framework.}
\newblock {\em Schizophrenia Bulletin}, 42(4):870--873.

\bibitem[Isvoranu et~al., 2016b]{isvoranu}
Isvoranu, A.~M., van Borkulo, C.~D., Boyette, L., Wigman, J. T.~W., Vinkers,
  C.~H., Borsboom, D., and {GROUP Investigators} (2016b).
\newblock {A Network Approach to Psychosis: Pathways between Childhood Trauma
  and Psychotic Symptoms}.
\newblock {\em Schizophrenia Bulletin}.
\newblock Advance Access published May 10, 2016.

\bibitem[Jacobucci et~al., 2016]{regsem}
Jacobucci, R., Grimm, K.~J., and McArdle, J.~J. (2016).
\newblock Regularized structural equation modeling.
\newblock {\em Structural Equation Modeling: A Multidisciplinary Journal},
  23(4):555--566.

\bibitem[J{\"o}reskog, 1967]{joreskog1967general}
J{\"o}reskog, K.~G. (1967).
\newblock A general approach to confirmatory maximum likelihood factor
  analysis.
\newblock {\em ETS Research Bulletin Series}, 1967(2):183--202.

\bibitem[Kaplan, 2000]{kaplan2008structural}
Kaplan, D. (2000).
\newblock {\em Structural equation modeling: Foundations and extensions}.
\newblock Sage, Thousand Oaks, CA, USA.

\bibitem[Kolaczyk, 2009]{kolaczyk2009statistical}
Kolaczyk, E.~D. (2009).
\newblock {\em Statistical analysis of network data}.
\newblock Springer, New York, NY, USA.

\bibitem[Koller and Friedman, 2009]{koller2009probabilistic}
Koller, D. and Friedman, N. (2009).
\newblock {\em Probabilistic graphical models: Principles and techniques}.
\newblock MIT press, Cambridge, MA, USA.

\bibitem[Kossakowski et~al., 2015]{kossakowski2015}
Kossakowski, J.~J., Epskamp, S., Kieffer, J.~M., van Borkulo, C.~D., Rhemtulla,
  M., and Borsboom, D. (2015).
\newblock The application of a network approach to health-related quality of
  life ({HRQoL}): Introducing a new method for assessing hrqol in healthy
  adults and cancer patient.
\newblock {\em Quality of Life Research}, 25:781--92.

\bibitem[Lauritzen, 1996]{lauritzen1996graphical}
Lauritzen, S.~L. (1996).
\newblock {\em Graphical models}.
\newblock Clarendon Press, Oxford, UK.

\bibitem[Lawley, 1940]{lawley1940vi}
Lawley, D.~N. (1940).
\newblock {VI}.---the estimation of factor loadings by the method of maximum
  likelihood.
\newblock {\em Proceedings of the Royal Society of Edinburgh}, 60(01):64--82.

\bibitem[Lord et~al., 1968]{lord1968statistical}
Lord, F.~M., Novick, M.~R., and Birnbaum, A. (1968).
\newblock {\em Statistical theories of mental test scores.}
\newblock Addison-Wesley, Oxford, UK.

\bibitem[MacCallum et~al., 1993]{maccallum1993problem}
MacCallum, R.~C., Wegener, D.~T., Uchino, B.~N., and Fabrigar, L.~R. (1993).
\newblock The problem of equivalent models in applications of covariance
  structure analysis.
\newblock {\em Psychological Bulletin}, 114(1):185--199.

\bibitem[Marsh et~al., 2014]{marsh2014exploratory}
Marsh, H.~W., Morin, A.~J., Parker, P.~D., and Kaur, G. (2014).
\newblock Exploratory structural equation modeling: An integration of the best
  features of exploratory and confirmatory factor analysis.
\newblock {\em Annual Review of Clinical Psychology}, 10:85--110.

\bibitem[Marsman et~al., 2015]{marsman2015bayesian}
Marsman, M., Maris, G., Bechger, T., and Glas, C. (2015).
\newblock Bayesian inference for low-rank ising networks.
\newblock {\em Scientific reports}, 5(9050):1--7.

\bibitem[McCrae and Costa, 1997]{mccrae1997personality}
McCrae, R.~R. and Costa, P.~T. (1997).
\newblock Personality trait structure as a human universal.
\newblock {\em American Psychologist}, 52(5):509--516.

\bibitem[McGill et~al., 1978]{mcgill1978variations}
McGill, R., Tukey, J.~W., and Larsen, W.~A. (1978).
\newblock Variations of box plots.
\newblock {\em The American Statistician}, 32(1):12--16.

\bibitem[McNally et~al., 2015]{mcnally2015mental}
McNally, R.~J., Robinaugh, D.~J., Wu, G.~W., Wang, L., Deserno, M.~K., and
  Borsboom, D. (2015).
\newblock Mental disorders as causal systems a network approach to
  posttraumatic stress disorder.
\newblock {\em Clinical Psychological Science}, 3(6):836--849.

\bibitem[Murphy, 2012]{murphy2012machine}
Murphy, K.~P. (2012).
\newblock {\em Machine learning: A probabilistic perspective}.
\newblock MIT press, Cambridge, MA, USA.

\bibitem[Neale et~al., 2016]{neale2015openmx}
Neale, M.~C., Hunter, M.~D., Pritikin, J.~N., Zahery, M., Brick, T.~R.,
  Kirkpatrick, R.~M., Estabrook, R., Bates, T.~C., Maes, H.~H., and Boker,
  S.~M. (2016).
\newblock Openmx 2.0: Extended structural equation and statistical modeling.
\newblock {\em Psychometrika}, 81(2):535--549.

\bibitem[Pearl, 2000]{pearl2000causality}
Pearl, J. (2000).
\newblock {\em Causality: Models, Reasoning, and Inference}.
\newblock Cambridge University Press, New York, NY.

\bibitem[Reckase, 2009]{reckase2009multidimensional}
Reckase, M.~D. (2009).
\newblock {\em Multidimensional item response theory}.
\newblock Springer, New York, NY, USA.

\bibitem[Revelle, 2010]{psych}
Revelle, W. (2010).
\newblock {\em {psych}: Procedures for Psychological, Psychometric, and
  Personality Research ({R} package version 1.0-93)}.
\newblock Northwestern University, Evanston, Illinois.

\bibitem[Rosa et~al., 2012]{rosa2012post}
Rosa, M., Friston, K., and Penny, W. (2012).
\newblock Post-hoc selection of dynamic causal models.
\newblock {\em Journal of neuroscience methods}, 208(1):66--78.

\bibitem[Schmittmann et~al., 2013]{schmittmann2013}
Schmittmann, V.~D., Cramer, A. O.~J., Waldorp, L.~J., Epskamp, S., Kievit,
  R.~A., and Borsboom, D. (2013).
\newblock {Deconstructing the construct: A network perspective on psychological
  phenomena}.
\newblock {\em New Ideas in Psychology}, 31(1):43--53.

\bibitem[Tibshirani, 1996]{tibshirani1996regression}
Tibshirani, R. (1996).
\newblock Regression shrinkage and selection via the lasso.
\newblock {\em Journal of the Royal Statistical Society. Series B
  (Methodological)}, 58:267--288.

\bibitem[van Borkulo et~al., 2014]{van2014new}
van Borkulo, C.~D., Borsboom, D., Epskamp, S., Blanken, T.~F., Boschloo, L.,
  Schoevers, R.~A., and Waldorp, L.~J. (2014).
\newblock A new method for constructing networks from binary data.
\newblock {\em Scientific Reports}, 4(5918):1--10.

\bibitem[van Borkulo et~al., 2015]{van2015association}
van Borkulo, C.~D., Boschloo, L., Borsboom, D., Penninx, B. W. J.~H., Waldorp,
  L.~J., and Schoevers, R.~A. (2015).
\newblock Association of symptom network structure with the course of
  depression.
\newblock {\em JAMA Psychiatry}, 72(12):1219--1226.

\bibitem[van~der Maas et~al., 2006]{van2006dynamical}
van~der Maas, H.~L., Dolan, C.~V., Grasman, R.~P., Wicherts, J.~M., Huizenga,
  H.~M., and Raijmakers, M.~E. (2006).
\newblock A dynamical model of general intelligence: The positive manifold of
  intelligence by mutualism.
\newblock {\em Psychological Review}, 113(4):842--861.

\bibitem[{Ware Jr} and Sherbourne, 1992]{ware1992}
{Ware Jr}, J.~E. and Sherbourne, C.~D. (1992).
\newblock {The MOS 36-item short-form health survey (SF-36): I. Conceptual
  framework and item selection}.
\newblock {\em Medical Care}, 30:473--483.

\bibitem[Wright, 1921]{wright1921correlation}
Wright, S. (1921).
\newblock Correlation and causation.
\newblock {\em Journal of agricultural research}, 20(7):557--585.

\bibitem[Wright, 1934]{wright1934method}
Wright, S. (1934).
\newblock The method of path coefficients.
\newblock {\em The Annals of Mathematical Statistics}, 5(3):161--215.

\bibitem[Yin and Li, 2011]{yin2011sparse}
Yin, J. and Li, H. (2011).
\newblock A sparse conditional gaussian graphical model for analysis of
  genetical genomics data.
\newblock {\em The Annals of Applied Statistics}, 5(4):2630--2650.

\bibitem[Zou and Hastie, 2005]{zou2005regularization}
Zou, H. and Hastie, T. (2005).
\newblock Regularization and variable selection via the elastic net.
\newblock {\em Journal of the Royal Statistical Society: Series B (Statistical
  Methodology)}, 67(2):301--320.

\bibitem[Zou et~al., 2007]{zou2007degrees}
Zou, H., Hastie, T., Tibshirani, R., et~al. (2007).
\newblock On the “degrees of freedom” of the lasso.
\newblock {\em The Annals of Statistics}, 35(5):2173--2192.

\end{thebibliography}

\end{document}